\documentclass[12pt,reqno]{article}

\usepackage[usenames]{color}
\usepackage{amssymb}
\usepackage{graphicx}
\usepackage{amscd}

\usepackage[colorlinks=true,
linkcolor=webgreen,
filecolor=webbrown,
citecolor=webgreen]{hyperref}

\definecolor{webgreen}{rgb}{0,.5,0}
\definecolor{webbrown}{rgb}{.6,0,0}

\usepackage{color}
\usepackage{float}

\usepackage{graphics,amsmath,amssymb}
\usepackage{amsthm}
\usepackage{amsfonts}
\usepackage{latexsym}

\begin{document}

\theoremstyle{plain}
\newtheorem{theorem}{Theorem}
\newtheorem{corollary}[theorem]{Corollary}
\newtheorem{lemma}[theorem]{Lemma}
\newtheorem{proposition}[theorem]{Proposition}
\newtheorem{obs}[theorem]{Observation}

\theoremstyle{definition}
\newtheorem{definition}[theorem]{Definition}
\newtheorem{example}[theorem]{Example}
\newtheorem{conjecture}[theorem]{Conjecture}

\theoremstyle{remark}
\newtheorem{remark}[theorem]{Remark}

\begin{center}
\vskip 1cm
{\LARGE\bf Representation of Dyck words in tensors that zipper merge contiguous integer compositions}

\vskip 1cm
\large
Italo J. Dejter\\
University of Puerto Rico\\
Rio Piedras, PR 00936-8377\\
\href{mailto:italo.dejter@gmail.com}{\tt italo.dejter@gmail.com} \\
\end{center}

\begin{abstract}
\noindent Let $0<k\in\mathbb{Z}$. We zipper merge integer compositions with sums $k$ and $k+1$, equal number of parts and initial entries equal at least to 1 and 2, respectively. This yields bitstrings with two initial zeros, $k-1$ remaining zeros and $k$ ones. Tensors whose entries are such bitstrings contain unique representations of all Dyck words of length $2k$. If rows and columns of such tensors are disposed in descending lexicographic order,
then their entries not representing Dyck words form disjoint unions of descending staircases corresponding to strict lower triangular submatrices.
\end{abstract}

\section{Zipper merging compositions of contiguous integers}\label{s1}

Let $0<k\in\mathbb{Z}$. An integer {\it composition} of $k$ is an ordered collection of one or more positive integers whose sum is $k$, \cite{Eger,HM1,HM2}. The summands of such sum are said to be the {\it parts} of the composition and the number of parts in a composition will be said to be its {\it length}. We will make use of the following counting facts.

\begin{theorem} \cite{HM2} The number of compositions of $k$ whose lengths are equal to $i$ is ${k-1\choose i-1}$, for $1\le i\le k$. The total number of compositions of $k$ is equal to $2^{k-1}$.
\end{theorem}

Let $P^k$ and $Q^k$ be the sets of compositions $(a_1,a_2,\ldots,a_i)$ of sum $k+1$ and $(b_1,b_2,\ldots, b_i)$ of sum $k$, respectively, where $a_1>1$. We subdivide $P^k$ and $Q^k$ into subsets $P^k_i$ and $Q^k_i$ according to the lengths $i$ of their participant compositions ($1\le i\le k$) and present these in descending lexicographic order, displaying both subsets jointly with the notation ${P^k_i\choose Q^k_i}$; see Table I, showing all the cases of ${P^k_i\choose Q^k_i}$ with $k=3,4,5$, where commas and parentheses of the compositions of the subsets $P^k_i$ and $Q^k_i$ are omitted.

For each case in Table I, we establish a corresponding associated matrix $P_i^k\otimes Q_i^k$, its rows and columns representing respective integer compositions $A=(a_1,a_2,\ldots, a_i)\in P^k_i$ and $B=(b_1,b_2,\ldots,b_i)\in Q^k_i$. Each such matrix is said to be the {\it tensor} $P_i^k\otimes Q_i^k$ that zipper merges the entries of the rows $A\in P_i^k$ with those of the columns $B\in Q_i^k$. The $(A,B)$-th entries of $P_i^k\otimes Q_i^k$ are zipper merging $2i$-tuples $A\otimes B=(a_1,b_1,a_2,b_2,\ldots,a_i,b_i)$. 

We replace the entries $a_j$ and $b_j$ of these $2i$-tuples $A\otimes B$ respectively by $a_j$ consecutive zeros followed by $b_j$ consecutive ones, for $j=1,\ldots,i$.
This yields bitstrings with two initial zeros, $k-1$ additional zeros and $k$ ones forming $(2k+1)$-tuples: 
 
 $$\mu_i^k(A,B)=\overbrace{00\cdots 0}^{a_1}\overbrace{11\cdots 1}^{b_1}\overbrace{00\cdots 0}^{a_2}\overbrace{11\cdots 1}^{b_2}\cdots\hspace*{5mm}\cdots\overbrace{00\cdots 0}^{a_i}\overbrace{11\cdots 1}^{b_i}=0|\nu_i^k(A,B).$$ 

{\center{\bf TABLE I}
$$\begin{array}{|l|l|l|l|l|}\hline
\!\!{P^3_1\choose Q^3_1}\!=\!{4\choose 3}\!\!&\!{P^3_2\choose Q^3_2}\!=\!{{31,22}\choose{21,12}}\!\!&\!\!{P^3_3\choose Q^3_3}\!=\!{{211}\choose{111}}\!\!&\!\!\!\!\!\!&\!\!\!\!\\\hline
\!\!{P^4_1\choose Q^4_1}\!=\!{5\choose 4}\!\!&\!\!{P^4_2\choose Q^4_2}={{41,32,23}\choose{31,22,13}}\!\!&\!\!{P^4_3\choose Q^4_3}\!=\!{{311,221,212}\choose{211,121,112}}\!\!&\!\!{P^4_3\choose Q^4_3}\!=\!{{2111}\choose{1111}}\!\!\!&\!\!\!\\\hline
\!\!{P^5_1\choose Q^5_1}\!=\!{6\choose 5}\!\!&\!\!{P^5_2\choose Q^5_2}\!=\!{{51,42,33,24}\choose{41,32,23,14}}\!\!&\!\!{P^5_3\choose Q^5_3}\!=\!{{411,321,312,231,222,123}\choose{311,221,212,131,122,113}}\!\!&\!\!{P^5_4\choose Q^5_4}\!=\!{{3111,2211,2121,2113}\choose{2111,1211,1121,1112}}\!\!&\!\!{P^5_5\choose Q^5_5}\!=\!{{21111}\choose{11111}}\!\\\hline\end{array}$$}

\begin{remark}\label{aqui} We consider bits $t_i^k(A,B)$ such that $t_i^k(A,B)=1$ if and only if $\nu_i^k(A,B)$ is a {\it Dyck word of length} $2k$, namely a $2k$-bitstring $0\cdots 1$ of weight $k$ such that in every prefix the number of 0-bits is at least the number of 1-bits \cite{Dejter0} (comparable to \cite{Hcs} \cite[p.~204]{Stanley}). Also, $t_i^k(A,B)=1$ if and only if $\mu_i^k(A,B)$ is equivalent to an {\it ordered $k$-edge tree} \cite[Remark 3]{Dejter0}, \cite{Dejter1}, \cite{gmn}, also known as ``plane tree with $k+1$ vertices" \cite[item~(e), p.~220]{Stanley}. 
\end{remark}

{\center{\bf TABLE II}
$$\begin{array}{|ccccc|cccccc|}\hline
&&T_1^3\!=\!\begin{bmatrix}1\end{bmatrix},\!&\!\!\!
T_2^3\!=\!\begin{bmatrix}1&1\\0 &1\end{bmatrix},\!&\!\!\!
T_3^3\!=\!\begin{bmatrix}1\end{bmatrix}\!&\!&\!
T_1^4\!=\!\begin{bmatrix}1\end{bmatrix},\!&\!\!
T_2^4\!=\!\begin{bmatrix}1&1&1\\0&1&1\\0&0&1\!\end{bmatrix}\!\!=
T_3^4,\!&\!
T_4^4\!=\!\begin{bmatrix}1\end{bmatrix}&&\!\\
\end{array}$$
\vspace*{-5mm}
$$\begin{array}{|ccccc|}\hline
T_1^5\!=\!\begin{bmatrix}1\end{bmatrix},\!&\!\!
T_2^5\!=\!\begin{bmatrix}1&1&1&1\!\\
0&1&1&1\!\\
0&0&1&1\!\\
0&0&0&1\!
\end{bmatrix},\!&\!\!
T_3^5\!=\!\begin{bmatrix}
1 &1&1&1&1&1\!\\
0  &1&1&1&1&\!1\!\!\\
0 & 0   &1& 0   &1&1\!\\
 0 &  0  & 0  &1&1&1\!\\
0& 0   & 0   & 0   &1&1\!\\
 0 & 0   & 0   &  0  & 0   &1\!
\end{bmatrix},\!&\!
T_4^5\!=\!\begin{bmatrix}1&1&1&1\!\\
0&1&1&1\!\\
0&0&1&1\!\\
0&0&0&1\!
\end{bmatrix},\!&\!\!
T_5^5\!=\!\begin{bmatrix}1\end{bmatrix}\!\!.\\\hline
\end{array}$$}

For $k>0$, $\mu_i^k(A,B)$ and $t_i^k(A,B)$ form ${k-1 \choose{i-1}}\times{k-1\choose{i-1}}$-tensors $U^k_1,\ldots,U^k_k$ and $T^k_1,\ldots,T^k_k$
associated to respective tensors $P^k_1\otimes Q^k_1,\ldots,P^k_k\otimes Q^k_k$. Let  $U^k_i=\{\mu^k_i(p,q)|p=1,\ldots,{k-1\choose i-1};q=1,\ldots,{k-1\choose i-1}\}$ and $T^k_i=\{t^k_i(p,q);p=1,\ldots,{k-1\choose i-1}|q=1,\ldots,{k-1\choose i-1}\}$. Denote the compositions assigned to the $p$-th row and to the $q$-th column of $T^k_i$ and $U^k_i$, both given in descending lexicographic order, by $A(p)=(a_1,\ldots,$ $a_i)=(a_1(p),\ldots,a_i(p))$ and $B(q)=(b_1,\ldots,b_i)$ $=(b_1(q),\ldots,$ $b_i(q))$, respectively. Table II shows those $T_i^k$ for the cases in Table I. Theorem~\ref{teo} and Corollary~\ref{theend} show that the entries of the tensors $T_i^k$ not representing Dyck words form disjoint unions of descending staircases corresponding to strict lower triangular submatrices.

{\center{\bf TABLE III}                   
$$\begin{array}{|llllllll|}\hline
(a_1^1)&=(4)&\in P^3_1&\mbox{ and }&(b_1^1)&=(3)&\in Q^3_1&\rightarrow 0000111\in U_1^3\\
(a_1^2,a_2^2)&=(3,1)&\in P^3_2&\mbox{ and }&(b_1^2,b_2^2)&=(2,1)&\in Q^3_2&\rightarrow 0001101\in U_2^3\\
(a_1^2,a_2^2)&=(3,1)&\in P^3_2&\mbox{ and }&(b_1^3,b_2^3)&=(1,2)&\in Q^3_2&\rightarrow 0001011\in U_2^3\\
(a_1^3,a_2^3)&=(2,2)&\in P^3_2&\mbox{ and }&(b_1^3,b_2^3)&=(1,2)&\in Q^3_3&\rightarrow 0010011\in U_2^3\\
(a_1^4,a_2^4,a_3^4)&=(2,1,1)&\in P_3^3&\mbox{ and }&(b_1^4,b_2^4,b_3^4)&=(1,1,1)&\in Q^3_3&\rightarrow 0010101\in U_3^3.\\\hline
\end{array}$$}

By force of Remark~\ref{aqui} above and Theorem~\ref{alla} below, we have: {\bf(1)} $T^k_1=[1]=T^k_k$; {\bf(2)} $t^k_i(p,q)=1$, if the partial sums $\sum_{j=1}^i[a_j(p)-b_j(q)]$ are all positive, and  $t^k_i(p,q)=0$, if at least one of those partial sums is not positive, for $1\le i\le k$; {\bf(3)} as a consequence, $t^k_i(p,q)=0$, for $1<i<k$ and for $1\le q<p\le k$, meaning $T^k_i$ is null below the main diagonal.

 \begin{theorem}\label{alla}
 {\bf(5)} The cases for which  
 $t^k_i(p,q)=1$ characterize the binary $(2k+1)$-tuples $\mu=\mu_i^k(A(p),B(q))$ representing ordered $k$-edge trees ${\mathcal T}(\mu)$, as well as their associated Dyck words $\nu_i^k(A(p),B(p))$ of length $2k$.\end{theorem}
 
 \begin{proof} The partial sums above, being all positive, imply the existence of a sequence of piecewise-lineal functions $\phi_j:[0,i]\rightarrow\mathbb{R}$ such that: {\bf(a)} $\phi_j(0)=0$ , ($j=0,\ldots,2k+1$), with $\phi_j$ positive on $(0,j]$, for $j\le 2k$; {\bf(b)} the function $\phi_{j+1}$ extends $\phi_j$ for $j=0,\ldots,2k$; and {\bf(c)} $\phi_{2k+1}(2k+1)=0$. This sequence of functions is obtained by replacing iteratively each non-starting entry 0, respectively 1, of $\mu$ by an ascending, respectively descending, unit-square diagonal that, in the $j$-th iteration, extends $\phi_j$ to the right into $\phi_{j+1}$ ($j\le 2k$). We also set right-to-left, bottom-to-top, increasing positive (lexical-matching \cite{KT}) integers in $(0,k]$, respectively asterisks, as labels on those ascending, respectively descending diagonals. As depicted for $k=3$ in the upper-right of Fig. 3,  the triangles and trapezoids resulting between contiguous integer heights in the graphic of $\phi_{2k+1}$ collapse horizontally onto the edges of ${\mathcal T}(\mu)$ looked-upon upside down, with: {\bf(i)} the root labeled 0 as the lowest vertex; {\bf(ii)} the remaining lexical colors numbering the highest vertices of corresponding collapsed ascending edges; and {\bf(iii)} asterisks denoting the collapsed descending edges. On each such ${\mathcal T}(\mu)$, depth preorder search done by going down an edge away from the root $r(\mu)$ of ${\mathcal T}(\mu)$, for each null entry, respectively going up an edge toward $r(\mu)$, for each unit entry, allows to reconstruct $\mu$ and so to recover $p$ and $q$ \cite{Dejter1,Dejter2,gmn,M,MN1,MN2}. On the other hand, \cite[ex.~6.19, p.~219]{Stanley} establishes that ordered $k$-edge trees, referred to as ``plane trees with" $k+1$ vertices in \cite[item~(e), p.~220]{Stanley}, are equivalent to the Dyck paths from $(0,0)$ to $(2k,0)$ of \cite[item~(i), p.221]{Stanley}, which in turn are equivalent to the Dick words of length $2k$ \cite{Dejter0}.
\end{proof}

{\center{\bf  TABLE IV}
 $$\begin{array}{|c|c|}\hline
 T_3^6=\begin{bmatrix}
1&1&1&1&1&1&1&1&1&1\\
0&1&1&1&1&1&1&1&1&1\\
0&0&1&0&1&1&0&1&1&1\\
0&0&0&1&1&1&1&1&1&1\\
0&0&0&0&1&1&0&1&1&1\\
0&0&0&0&0&1&0&0&1&1\\
0&0&0&0&0&0&1&1&1&1\\
0&0&0&0&0&0&0&1&1&1\\
0&0&0&0&0&0&0&0&1&1\\
0&0&0&0&0&0&0&0&0&1\\
 \end{bmatrix}&
 T_4^6=\begin{bmatrix}
1&1&1&1&1&1&1&1&1&1\\
0&1&1&1&1&1&1&1&1&1\\
0&0&1&1&0&1&1&1&1&1\\
0&0&0&1&0&0&1&0&1&1\\
0&0&0&0&1&1&1&1&1&1\\
0&0&0&0&0&1&1&1&1&1\\
0&0&0&0&0&0&1&0&1&1\\
0&0&0&0&0&0&0&1&1&1\\
0&0&0&0&0&0&0&0&1&1\\
0&0&0&0&0&0&0&0&0&1\\
 \end{bmatrix}\\\hline
 \end{array}$$}

By replacing the compositions associated to the rows and columns of the three tensors $T_i^k$ shown on the upper-left of Table II and zipper merging their information as suggested above in place of their nonzero entries, we get tensors $U_i^k$ whose entries are the $\mu_i^k(A(p),B(q))$:
\begin{eqnarray}\begin{array}{lll}
T_1^3\rightarrow\!\begin{bmatrix}^{\;\;|3}_{4|0000111}\end{bmatrix}=U_1^3,\!&\!
T_2^3\rightarrow\!\begin{bmatrix}^{\;\;\;\;|21}_{31|0001101}&^{|12}_{|0001011}\\
^{22|----.}&^{|0010011}
\end{bmatrix}=U_2^3,&
T_3^3\rightarrow\!\begin{bmatrix}^{\;\;\;\;\;\,|111}_{211|0010101}\end{bmatrix}=U_3^3.
\end{array}\end{eqnarray}
Those compositions are zipper merged in Table III (via arrows ``$\rightarrow$'') into the binary $(2k+1)$-tuples in display (1) presented (via symbols ``$\in$'') as belonging to their tensors $U_i^k$, representing the ordered trees suggested in the upper-right of Fig.~3. 

{\center{\bf TABLE V}
$$\begin{array}{|ccccl|l|lcccccc|}\hline
 &&&& T_3^5\!=\!\begin{bmatrix}
^{\bullet\,\,\bullet\,\,\bullet\,\,\bullet\,\,\bullet\,\,\bullet}_
{\circ\,\,\bullet\,\,\bullet\,\,\bullet\,\,\bullet\,\,\bullet}\\
^{\circ\,\,\circ\,\,\bullet\,\,\circ\,\,\bullet\,\,\bullet}_
{\circ\,\,\circ\,\,\circ\,\,\bullet\,\,\bullet\,\,\bullet}\\
^{\circ\,\,\circ\,\,\circ\,\,\circ\,\,\bullet\,\,\bullet}_
{\circ\,\,\circ\,\,\circ\,\,\circ\,\,\circ\,\,\bullet}
\end{bmatrix}&\,
 T_3^6=\begin{bmatrix}
^{\,\,\bullet\,\,\bullet\,\,\bullet\,\,\bullet\,\,\bullet\,\,\bullet\,\,\bullet\,\,\bullet\,\,\bullet\,\,\bullet\,\,}
_{\,\,\circ\,\,\bullet\,\,\bullet\,\,\bullet\,\,\bullet\,\,\bullet\,\,\bullet\,\,\bullet\,\,\bullet\,\,\bullet\,\,}\\
^{\,\,\circ\,\,\circ\,\,\bullet\,\,\circ\,\,\bullet\,\,\bullet\,\,\circ\,\,\bullet\,\,\bullet\,\,\bullet\,\,}
_{\,\,\circ\,\,\circ\,\,\circ\,\,\bullet\,\,\bullet\,\,\bullet\,\,\bullet\,\,\bullet\,\,\bullet\,\,\bullet\,\,}\\
^{\,\,\circ\,\,\circ\,\,\circ\,\,\circ\,\,\bullet\,\,\bullet\,\,\circ\,\,\bullet\,\,\bullet\,\,\bullet\,\,}
_{\,\,\circ\,\,\circ\,\,\circ\,\,\circ\,\,\circ\,\,\bullet\,\,\circ\,\,\circ\,\,\bullet\,\,\bullet\,\,}\\
^{\,\,\circ\,\,\circ\,\,\circ\,\,\circ\,\,\circ\,\,\circ\,\,\bullet\,\,\bullet\,\,\bullet\,\,\bullet\,\,}
_{\,\,\circ\,\,\circ\,\,\circ\,\,\circ\,\,\circ\,\,\circ\,\,\circ\,\,\bullet\,\,\bullet\,\,\bullet\,\,}\\
^{\,\,\circ\,\,\circ\,\,\circ\,\,\circ\,\,\circ\,\,\circ\,\,\circ\,\,\circ\,\,\bullet\,\,\bullet\,\,}
_{\,\,\circ\,\,\circ\,\,\circ\,\,\circ\,\,\circ\,\,\circ\,\,\circ\,\,\circ\,\,\circ\,\,\bullet\,\,}\\
 \end{bmatrix}&
 T_4^6=\begin{bmatrix}
^{\,\,\bullet\,\,\bullet\,\,\bullet\,\,\bullet\,\,\bullet\,\,\bullet\,\,\bullet\,\,\bullet\,\,\bullet\,\,\bullet\,\,}
_{\,\,\circ\,\,\bullet\,\,\bullet\,\,\bullet\,\,\bullet\,\,\bullet\,\,\bullet\,\,\bullet\,\,\bullet\,\,\bullet\,\,}\\
^{\,\,\circ\,\,\circ\,\,\bullet\,\,\bullet\,\,\circ\,\,\bullet\,\,\bullet\,\,\bullet\,\,\bullet\,\,\bullet\,\,}
_{\,\,\circ\,\,\circ\,\,\circ\,\,\bullet\,\,\circ\,\,\circ\,\,\bullet\,\,\circ\,\,\bullet\,\,\bullet\,\,}\\
^{\,\,\circ\,\,\circ\,\,\circ\,\,\circ\,\,\bullet\,\,\bullet\,\,\bullet\,\,\bullet\,\,\bullet\,\,\bullet\,\,}
_{\,\,\circ\,\,\circ\,\,\circ\,\,\circ\,\,\circ\,\,\bullet\,\,\bullet\,\,\bullet\,\,\bullet\,\,\bullet\,\,}\\
^{\,\,\circ\,\,\circ\,\,\circ\,\,\circ\,\,\circ\,\,\circ\,\,\bullet\,\,\circ\,\,\bullet\,\,\bullet\,\,}
_{\,\,\circ\,\,\circ\,\,\circ\,\,\circ\,\,\circ\,\,\circ\,\,\circ\,\,\bullet\,\,\bullet\,\,\bullet\,\,}\\
^{\,\,\circ\,\,\circ\,\,\circ\,\,\circ\,\,\circ\,\,\circ\,\,\circ\,\,\circ\,\,\bullet\,\,\bullet\,\,}
_{\,\,\circ\,\,\circ\,\,\circ\,\,\circ\,\,\circ\,\,\circ\,\,\circ\,\,\circ\,\,\circ\,\,\bullet\,\,}\\
 \end{bmatrix}&&&&&&\\
 \end{array}$$
 \vspace*{-3mm}
$$\begin{array}{|l|l|llll|}\hline
T_3^7=&T_4^7=&T_5^7=&&&\\
\begin{bmatrix}
^{\bullet\,\,\bullet\,\,\bullet\,\,\bullet\,\,\bullet\,\,\bullet\,\,\bullet\,\,\bullet\,\,\bullet\,\,\bullet\,\,\bullet\,\,\bullet\,\,\bullet\,\,\bullet\,\,\bullet\,\,}
_{\circ\,\,\bullet\,\,\bullet\,\,\bullet\,\,\bullet\,\,\bullet\,\,\bullet\,\,\bullet\,\,\bullet\,\,\bullet\,\,\bullet\,\,\bullet\,\,\bullet\,\,\bullet\,\,\bullet\,\,}\\
^{\circ\,\,\circ\,\,\bullet\,\,\circ\,\,\bullet\,\,\bullet\,\,\circ\,\,\bullet\,\,\bullet\,\,\bullet\,\,\circ\,\,\bullet\,\,\bullet\,\,\bullet\,\,\bullet\,\,}
_{\circ\,\,\circ\,\,\circ\,\,\bullet\,\,\bullet\,\,\bullet\,\,\bullet\,\,\bullet\,\,\bullet\,\,\bullet\,\,\bullet\,\,\bullet\,\,\bullet\,\,\bullet\,\,\bullet\,\,}\\
^{\circ\,\,\circ\,\,\circ\,\,\circ\,\,\bullet\,\,\bullet\,\,\circ\,\,\bullet\,\,\bullet\,\,\bullet\,\,\circ\,\,\bullet\,\,\bullet\,\,\bullet\,\,\bullet\,\,}
_{\circ\,\,\circ\,\,\circ\,\,\circ\,\,\circ\,\,\bullet\,\,\circ\,\,\circ\,\,\bullet\,\,\bullet\,\,\circ\,\,\circ\,\,\bullet\,\,\bullet\,\,\bullet\,\,}\\
^{\circ\,\,\circ\,\,\circ\,\,\circ\,\,\circ\,\,\circ\,\,\bullet\,\,\bullet\,\,\bullet\,\,\bullet\,\,\bullet\,\,\bullet\,\,\bullet\,\,\bullet\,\,\bullet\,\,}
_{\circ\,\,\circ\,\,\circ\,\,\circ\,\,\circ\,\,\circ\,\,\circ\,\,\bullet\,\,\bullet\,\,\bullet\,\,\circ\,\,\bullet\,\,\bullet\,\,\bullet\,\,\bullet\,\,}\\
^{\circ\,\,\circ\,\,\circ\,\,\circ\,\,\circ\,\,\circ\,\,\circ\,\,\circ\,\,\bullet\,\,\bullet\,\,\circ\,\,\circ\,\,\bullet\,\,\bullet\,\,\bullet\,\,}
_{\circ\,\,\circ\,\,\circ\,\,\circ\,\,\circ\,\,\circ\,\,\circ\,\,\circ\,\,\circ\,\,\bullet\,\,\circ\,\,\circ\,\,\circ\,\,\bullet\,\,\bullet\,\,}\\
^{\circ\,\,\circ\,\,\circ\,\,\circ\,\,\circ\,\,\circ\,\,\circ\,\,\circ\,\,\circ\,\,\circ\,\,\bullet\,\,\bullet\,\,\bullet\,\,\bullet\,\,\bullet\,\,}
_{\circ\,\,\circ\,\,\circ\,\,\circ\,\,\circ\,\,\circ\,\,\circ\,\,\circ\,\,\circ\,\,\circ\,\,\circ\,\,\bullet\,\,\bullet\,\,\bullet\,\,\bullet\,\,}\\
^{\circ\,\,\circ\,\,\circ\,\,\circ\,\,\circ\,\,\circ\,\,\circ\,\,\circ\,\,\circ\,\,\circ\,\,\circ\,\,\circ\,\,\bullet\,\,\bullet\,\,\bullet\,\,}
_{\circ\,\,\circ\,\,\circ\,\,\circ\,\,\circ\,\,\circ\,\,\circ\,\,\circ\,\,\circ\,\,\circ\,\,\circ\,\,\circ\,\,\circ\,\,\bullet\,\,\bullet\,\,}\\
^{\circ\,\,\circ\,\,\circ\,\,\circ\,\,\circ\,\,\circ\,\,\circ\,\,\circ\,\,\circ\,\,\circ\,\,\circ\,\,\circ\,\,\circ\,\,\circ\,\,\bullet\,\,}\\\end{bmatrix}&
\begin{bmatrix}
^{\bullet\,\,\bullet\,\,\bullet\,\,\bullet\,\,\bullet\,\,\bullet\,\,\bullet\,\,\bullet\,\,\bullet\,\,\bullet\,\,\bullet\,\,\bullet\,\,\bullet\,\,\bullet\,\,\bullet\,\,\bullet\,\,\bullet\,\,\bullet\,\,\bullet\,\,\bullet\,\,}
_{\circ\,\,\bullet\,\,\bullet\,\,\bullet\,\,\bullet\,\,\bullet\,\,\bullet\,\,\bullet\,\,\bullet\,\,\bullet\,\,\bullet\,\,\bullet\,\,\bullet\,\,\bullet\,\,\bullet\,\,\bullet\,\,\bullet\,\,\bullet\,\,\bullet\,\,\bullet\,\,}\\
^{\circ\,\,\circ\,\,\bullet\,\,\bullet\,\,\circ\,\,\bullet\,\,\bullet\,\,\bullet\,\,\bullet\,\,\bullet\,\,\circ\,\,\bullet\,\,\bullet\,\,\bullet\,\,\bullet\,\,\bullet\,\,\bullet\,\,\bullet\,\,\bullet\,\,\bullet\,\,}
_{\circ\,\,\circ\,\,\circ\,\,\bullet\,\,\circ\,\,\circ\,\,\bullet\,\,\circ\,\,\bullet\,\,\bullet\,\,\circ\,\,\circ\,\,\bullet\,\,\circ\,\,\bullet\,\,\bullet\,\,\circ\,\,\bullet\,\,\bullet\,\,\bullet\,\,}\\
^{\circ\,\,\circ\,\,\circ\,\,\circ\,\,\bullet\,\,\bullet\,\,\bullet\,\,\bullet\,\,\bullet\,\,\bullet\,\,\bullet\,\,\bullet\,\,\bullet\,\,\bullet\,\,\bullet\,\,\bullet\,\,\bullet\,\,\bullet\,\,\bullet\,\,\bullet\,\,}
_{\circ\,\,\circ\,\,\circ\,\,\circ\,\,\circ\,\,\bullet\,\,\bullet\,\,\bullet\,\,\bullet\,\,\bullet\,\,\circ\,\,\bullet\,\,\bullet\,\,\bullet\,\,\bullet\,\,\bullet\,\,\bullet\,\,\bullet\,\,\bullet\,\,\bullet\,\,}\\
^{\circ\,\,\circ\,\,\circ\,\,\circ\,\,\circ\,\,\circ\,\,\bullet\,\,\circ\,\,\bullet\,\,\bullet\,\,\circ\,\,\circ\,\,\bullet\,\,\circ\,\,\bullet\,\,\bullet\,\,\circ\,\,\bullet\,\,\bullet\,\,\bullet\,\,}
_{\circ\,\,\circ\,\,\circ\,\,\circ\,\,\circ\,\,\circ\,\,\circ\,\,\bullet\,\,\bullet\,\,\bullet\,\,\circ\,\,\circ\,\,\circ\,\,\bullet\,\,\bullet\,\,\bullet\,\,\bullet\,\,\bullet\,\,\bullet\,\,\bullet\,\,}\\
^{\circ\,\,\circ\,\,\circ\,\,\circ\,\,\circ\,\,\circ\,\,\circ\,\,\circ\,\,\bullet\,\,\bullet\,\,\circ\,\,\circ\,\,\circ\,\,\circ\,\,\bullet\,\,\bullet\,\,\circ\,\,\bullet\,\,\bullet\,\,\bullet\,\,}
_{\circ\,\,\circ\,\,\circ\,\,\circ\,\,\circ\,\,\circ\,\,\circ\,\,\circ\,\,\circ\,\,\bullet\,\,\circ\,\,\circ\,\,\circ\,\,\circ\,\,\circ\,\,\bullet\,\,\circ\,\,\circ\,\,\bullet\,\,\bullet\,\,}\\
^{\circ\,\,\circ\,\,\circ\,\,\circ\,\,\circ\,\,\circ\,\,\circ\,\,\circ\,\,\circ\,\,\circ\,\,\bullet\,\,\bullet\,\,\bullet\,\,\bullet\,\,\bullet\,\,\bullet\,\,\bullet\,\,\bullet\,\,\bullet\,\,\bullet\,\,}
_{\circ\,\,\circ\,\,\circ\,\,\circ\,\,\circ\,\,\circ\,\,\circ\,\,\circ\,\,\circ\,\,\circ\,\,\circ\,\,\bullet\,\,\bullet\,\,\bullet\,\,\bullet\,\,\bullet\,\,\bullet\,\,\bullet\,\,\bullet\,\,\bullet\,\,}\\
^{\circ\,\,\circ\,\,\circ\,\,\circ\,\,\circ\,\,\circ\,\,\circ\,\,\circ\,\,\circ\,\,\circ\,\,\circ\,\,\circ\,\,\bullet\,\,\circ\,\,\bullet\,\,\bullet\,\,\circ\,\,\bullet\,\,\bullet\,\,\bullet\,\,}
_{\circ\,\,\circ\,\,\circ\,\,\circ\,\,\circ\,\,\circ\,\,\circ\,\,\circ\,\,\circ\,\,\circ\,\,\circ\,\,\circ\,\,\circ\,\,\bullet\,\,\bullet\,\,\bullet\,\,\bullet\,\,\bullet\,\,\bullet\,\,\bullet\,\,}\\
^{\circ\,\,\circ\,\,\circ\,\,\circ\,\,\circ\,\,\circ\,\,\circ\,\,\circ\,\,\circ\,\,\circ\,\,\circ\,\,\circ\,\,\circ\,\,\circ\,\,\bullet\,\,\bullet\,\,\circ\,\,\bullet\,\,\bullet\,\,\bullet\,\,}
_{\circ\,\,\circ\,\,\circ\,\,\circ\,\,\circ\,\,\circ\,\,\circ\,\,\circ\,\,\circ\,\,\circ\,\,\circ\,\,\circ\,\,\circ\,\,\circ\,\,\circ\,\,\bullet\,\,\circ\,\,\circ\,\,\bullet\,\,\bullet\,\,}\\
^{\circ\,\,\circ\,\,\circ\,\,\circ\,\,\circ\,\,\circ\,\,\circ\,\,\circ\,\,\circ\,\,\circ\,\,\circ\,\,\circ\,\,\circ\,\,\circ\,\,\circ\,\,\circ\,\,\bullet\,\,\bullet\,\,\bullet\,\,\bullet\,\,}
_{\circ\,\,\circ\,\,\circ\,\,\circ\,\,\circ\,\,\circ\,\,\circ\,\,\circ\,\,\circ\,\,\circ\,\,\circ\,\,\circ\,\,\circ\,\,\circ\,\,\circ\,\,\circ\,\,\circ\,\,\bullet\,\,\bullet\,\,\bullet\,\,}\\
^{\circ\,\,\circ\,\,\circ\,\,\circ\,\,\circ\,\,\circ\,\,\circ\,\,\circ\,\,\circ\,\,\circ\,\,\circ\,\,\circ\,\,\circ\,\,\circ\,\,\circ\,\,\circ\,\,\circ\,\,\circ\,\,\bullet\,\,\bullet\,\,}
_{\circ\,\,\circ\,\,\circ\,\,\circ\,\,\circ\,\,\circ\,\,\circ\,\,\circ\,\,\circ\,\,\circ\,\,\circ\,\,\circ\,\,\circ\,\,\circ\,\,\circ\,\,\circ\,\,\circ\,\,\circ\,\,\circ\,\,\bullet\,\,}\\
\end{bmatrix}&
\begin{bmatrix}
^{\bullet\,\,\bullet\,\,\bullet\,\,\bullet\,\,\bullet\,\,\bullet\,\,\bullet\,\,\bullet\,\,\bullet\,\,\bullet\,\,\bullet\,\,\bullet\,\,\bullet\,\,\bullet\,\,\bullet\,\,}
_{\circ\,\,\bullet\,\,\bullet\,\,\bullet\,\,\bullet\,\,\bullet\,\,\bullet\,\,\bullet\,\,\bullet\,\,\bullet\,\,\bullet\,\,\bullet\,\,\bullet\,\,\bullet\,\,\bullet\,\,}\\
^{\circ\,\,\circ\,\,\bullet\,\,\bullet\,\,\bullet\,\,\circ\,\,\bullet\,\,\bullet\,\,\bullet\,\,\bullet\,\,\bullet\,\,\bullet\,\,\bullet\,\,\bullet\,\,\bullet\,\,}
_{\circ\,\,\circ\,\,\circ\,\,\bullet\,\,\bullet\,\,\circ\,\,\circ\,\,\bullet\,\,\bullet\,\,\circ\,\,\bullet\,\,\bullet\,\,\bullet\,\,\bullet\,\,\bullet\,\,}\\
^{\circ\,\,\circ\,\,\circ\,\,\circ\,\,\bullet\,\,\circ\,\,\circ\,\,\circ\,\,\bullet\,\,\circ\,\,\circ\,\,\bullet\,\,\circ\,\,\bullet\,\,\bullet\,\,}
_{\circ\,\,\circ\,\,\circ\,\,\circ\,\,\circ\,\,\bullet\,\,\bullet\,\,\bullet\,\,\bullet\,\,\bullet\,\,\bullet\,\,\bullet\,\,\bullet\,\,\bullet\,\,\bullet\,\,}\\
^{\circ\,\,\circ\,\,\circ\,\,\circ\,\,\circ\,\,\circ\,\,\bullet\,\,\bullet\,\,\bullet\,\,\bullet\,\,\bullet\,\,\bullet\,\,\bullet\,\,\bullet\,\,\bullet\,\,}
_{\circ\,\,\circ\,\,\circ\,\,\circ\,\,\circ\,\,\circ\,\,\circ\,\,\bullet\,\,\bullet\,\,\circ\,\,\bullet\,\,\bullet\,\,\bullet\,\,\bullet\,\,\bullet\,\,}\\
^{\circ\,\,\circ\,\,\circ\,\,\circ\,\,\circ\,\,\circ\,\,\circ\,\,\circ\,\,\bullet\,\,\circ\,\,\circ\,\,\bullet\,\,\circ\,\,\bullet\,\,\bullet\,\,}
_{\circ\,\,\circ\,\,\circ\,\,\circ\,\,\circ\,\,\circ\,\,\circ\,\,\circ\,\,\circ\,\,\bullet\,\,\bullet\,\,\bullet\,\,\bullet\,\,\bullet\,\,\bullet\,\,}\\
^{\circ\,\,\circ\,\,\circ\,\,\circ\,\,\circ\,\,\circ\,\,\circ\,\,\circ\,\,\circ\,\,\circ\,\,\bullet\,\,\bullet\,\,\bullet\,\,\bullet\,\,\bullet\,\,}
_{\circ\,\,\circ\,\,\circ\,\,\circ\,\,\circ\,\,\circ\,\,\circ\,\,\circ\,\,\circ\,\,\circ\,\,\circ\,\,\bullet\,\,\circ\,\,\bullet\,\,\bullet\,\,}\\
^{\circ\,\,\circ\,\,\circ\,\,\circ\,\,\circ\,\,\circ\,\,\circ\,\,\circ\,\,\circ\,\,\circ\,\,\circ\,\,\circ\,\,\bullet\,\,\bullet\,\,\bullet\,\,}
_{\circ\,\,\circ\,\,\circ\,\,\circ\,\,\circ\,\,\circ\,\,\circ\,\,\circ\,\,\circ\,\,\circ\,\,\circ\,\,\circ\,\,\circ\,\,\bullet\,\,\bullet\,\,}\\
^{\circ\,\,\circ\,\,\circ\,\,\circ\,\,\circ\,\,\circ\,\,\circ\,\,\circ\,\,\circ\,\,\circ\,\,\circ\,\,\circ\,\,\circ\,\,\circ\,\,\bullet\,\,}\\\end{bmatrix}&&&
\end{array}$$
\vspace*{-2mm}
$$\begin{array}{|l|lc|}\hline
T_3^8=&T_4^8=&\\
\begin{bmatrix}
^{\bullet\,\,\bullet\,\,\bullet\,\,\bullet\,\,\bullet\,\,\bullet\,\,\bullet\,\,\bullet\,\,\bullet\,\,\bullet\,\,\bullet\,\,\bullet\,\,\bullet\,\,\bullet\,\,\bullet\,\,\bullet\,\,\bullet\,\,\bullet\,\,\bullet\,\,\bullet\,\,\bullet\,\,}_
{\circ\,\,\bullet\,\,\bullet\,\,\bullet\,\,\bullet\,\,\bullet\,\,\bullet\,\,\bullet\,\,\bullet\,\,\bullet\,\,\bullet\,\,\bullet\,\,\bullet\,\,\bullet\,\,\bullet\,\,\bullet\,\,\bullet\,\,\bullet\,\,\bullet\,\,\bullet\,\,\bullet\,\,}\\
^{\circ\,\,\circ\,\,\bullet\,\,\circ\,\,\bullet\,\,\bullet\,\,\circ\,\,\bullet\,\,\bullet\,\,\bullet\,\,\circ\,\,\bullet\,\,\bullet\,\,\bullet\,\,\bullet\,\,\circ\,\,\bullet\,\,\bullet\,\,\bullet\,\,\bullet\,\,\bullet\,\,}_
{\circ\,\,\circ\,\,\circ\,\,\bullet\,\,\bullet\,\,\bullet\,\,\bullet\,\,\bullet\,\,\bullet\,\,\bullet\,\,\bullet\,\,\bullet\,\,\bullet\,\,\bullet\,\,\bullet\,\,\bullet\,\,\bullet\,\,\bullet\,\,\bullet\,\,\bullet\,\,\bullet\,\,}\\
^{\circ\,\,\circ\,\,\circ\,\,\circ\,\,\bullet\,\,\bullet\,\,\circ\,\,\bullet\,\,\bullet\,\,\bullet\,\,\circ\,\,\bullet\,\,\bullet\,\,\bullet\,\,\bullet\,\,\circ\,\,\bullet\,\,\bullet\,\,\bullet\,\,\bullet\,\,\bullet\,\,}_
{\circ\,\,\circ\,\,\circ\,\,\circ\,\,\circ\,\,\bullet\,\,\circ\,\,\circ\,\,\bullet\,\,\bullet\,\,\circ\,\,\circ\,\,\bullet\,\,\bullet\,\,\bullet\,\,\circ\,\,\circ\,\,\bullet\,\,\bullet\,\,\bullet\,\,\bullet\,\,}\\
^{\circ\,\,\circ\,\,\circ\,\,\circ\,\,\circ\,\,\circ\,\,\bullet\,\,\bullet\,\,\bullet\,\,\bullet\,\,\bullet\,\,\bullet\,\,\bullet\,\,\bullet\,\,\bullet\,\,\bullet\,\,\bullet\,\,\bullet\,\,\bullet\,\,\bullet\,\,\bullet\,\,}_
{\circ\,\,\circ\,\,\circ\,\,\circ\,\,\circ\,\,\circ\,\,\circ\,\,\bullet\,\,\bullet\,\,\bullet\,\,\circ\,\,\bullet\,\,\bullet\,\,\bullet\,\,\bullet\,\,\circ\,\,\bullet\,\,\bullet\,\,\bullet\,\,\bullet\,\,\bullet\,\,}\\
^{\circ\,\,\circ\,\,\circ\,\,\circ\,\,\circ\,\,\circ\,\,\circ\,\,\circ\,\,\bullet\,\,\bullet\,\,\circ\,\,\circ\,\,\bullet\,\,\bullet\,\,\bullet\,\,\circ\,\,\circ\,\,\bullet\,\,\bullet\,\,\bullet\,\,\bullet\,\,}_
{\circ\,\,\circ\,\,\circ\,\,\circ\,\,\circ\,\,\circ\,\,\circ\,\,\circ\,\,\circ\,\,\bullet\,\,\circ\,\,\circ\,\,\circ\,\,\bullet\,\,\bullet\,\,\circ\,\,\circ\,\,\circ\,\,\bullet\,\,\bullet\,\,\bullet\,\,}\\
^{\circ\,\,\circ\,\,\circ\,\,\circ\,\,\circ\,\,\circ\,\,\circ\,\,\circ\,\,\circ\,\,\circ\,\,\bullet\,\,\bullet\,\,\bullet\,\,\bullet\,\,\bullet\,\,\bullet\,\,\bullet\,\,\bullet\,\,\bullet\,\,\bullet\,\,\bullet\,\,}
_{\circ\,\,\circ\,\,\circ\,\,\circ\,\,\circ\,\,\circ\,\,\circ\,\,\circ\,\,\circ\,\,\circ\,\,\circ\,\,\bullet\,\,\bullet\,\,\bullet\,\,\bullet\,\,\circ\,\,\bullet\,\,\bullet\,\,\bullet\,\,\bullet\,\,\bullet\,\,}\\
^{\circ\,\,\circ\,\,\circ\,\,\circ\,\,\circ\,\,\circ\,\,\circ\,\,\circ\,\,\circ\,\,\circ\,\,\circ\,\,\circ\,\,\bullet\,\,\bullet\,\,\bullet\,\,\circ\,\,\circ\,\,\bullet\,\,\bullet\,\,\bullet\,\,\bullet\,\,}_
{\circ\,\,\circ\,\,\circ\,\,\circ\,\,\circ\,\,\circ\,\,\circ\,\,\circ\,\,\circ\,\,\circ\,\,\circ\,\,\circ\,\,\circ\,\,\bullet\,\,\bullet\,\,\circ\,\,\circ\,\,\circ\,\,\bullet\,\,\bullet\,\,\bullet\,\,}\\
^{\circ\,\,\circ\,\,\circ\,\,\circ\,\,\circ\,\,\circ\,\,\circ\,\,\circ\,\,\circ\,\,\circ\,\,\circ\,\,\circ\,\,\circ\,\,\circ\,\,\bullet\,\,\circ\,\,\circ\,\,\circ\,\,\circ\,\,\bullet\,\,\bullet\,\,}_
{\circ\,\,\circ\,\,\circ\,\,\circ\,\,\circ\,\,\circ\,\,\circ\,\,\circ\,\,\circ\,\,\circ\,\,\circ\,\,\circ\,\,\circ\,\,\circ\,\,\circ\,\,\bullet\,\,\bullet\,\,\bullet\,\,\bullet\,\,\bullet\,\,\bullet\,\,}\\
^{\circ\,\,\circ\,\,\circ\,\,\circ\,\,\circ\,\,\circ\,\,\circ\,\,\circ\,\,\circ\,\,\circ\,\,\circ\,\,\circ\,\,\circ\,\,\circ\,\,\circ\,\,\circ\,\,\bullet\,\,\bullet\,\,\bullet\,\,\bullet\,\,\bullet\,\,}_
{\circ\,\,\circ\,\,\circ\,\,\circ\,\,\circ\,\,\circ\,\,\circ\,\,\circ\,\,\circ\,\,\circ\,\,\circ\,\,\circ\,\,\circ\,\,\circ\,\,\circ\,\,\circ\,\,\circ\,\,\bullet\,\,\bullet\,\,\bullet\,\,\bullet\,\,}\\
^{\circ\,\,\circ\,\,\circ\,\,\circ\,\,\circ\,\,\circ\,\,\circ\,\,\circ\,\,\circ\,\,\circ\,\,\circ\,\,\circ\,\,\circ\,\,\circ\,\,\circ\,\,\circ\,\,\circ\,\,\circ\,\,\bullet\,\,\bullet\,\,\bullet\,\,}_
{\circ\,\,\circ\,\,\circ\,\,\circ\,\,\circ\,\,\circ\,\,\circ\,\,\circ\,\,\circ\,\,\circ\,\,\circ\,\,\circ\,\,\circ\,\,\circ\,\,\circ\,\,\circ\,\,\circ\,\,\circ\,\,\circ\,\,\bullet\,\,\bullet\,\,}\\
^{\circ\,\,\circ\,\,\circ\,\,\circ\,\,\circ\,\,\circ\,\,\circ\,\,\circ\,\,\circ\,\,\circ\,\,\circ\,\,\circ\,\,\circ\,\,\circ\,\,\circ\,\,\circ\,\,\circ\,\,\circ\,\,\circ\,\,\circ\,\,\bullet\,\,}
\end{bmatrix}&\,\;
\begin{bmatrix}
^{\bullet\,\,\bullet\,\,\bullet\,\,\bullet\,\,\bullet\,\,\bullet\,\,\bullet\,\,\bullet\,\,\bullet\,\,\bullet\,\,\bullet\,\,\bullet\,\,\bullet\,\,\bullet\,\,\bullet\,\,\bullet\,\,\bullet\,\,\bullet\,\,\bullet\,\,\bullet\,\,\bullet\,\,\bullet\,\,\bullet\,\,\bullet\,\,\bullet\,\,\bullet\,\,\bullet\,\,\bullet\,\,\bullet\,\,\bullet\,\,\bullet\,\,\bullet\,\,\bullet\,\,\bullet\,\,\bullet\,\,}
_{\circ\,\,\bullet\,\,\bullet\,\,\bullet\,\,\bullet\,\,\bullet\,\,\bullet\,\,\bullet\,\,\bullet\,\,\bullet\,\,\bullet\,\,\bullet\,\,\bullet\,\,\bullet\,\,\bullet\,\,\bullet\,\,\bullet\,\,\bullet\,\,\bullet\,\,\bullet\,\,\bullet\,\,\bullet\,\,\bullet\,\,\bullet\,\,\bullet\,\,\bullet\,\,\bullet\,\,\bullet\,\,\bullet\,\,\bullet\,\,\bullet\,\,\bullet\,\,\bullet\,\,\bullet\,\,\bullet\,\,}\\
^{\circ\,\,\circ\,\,\bullet\,\,\bullet\,\,\circ\,\,\bullet\,\,\bullet\,\,\bullet\,\,\bullet\,\,\bullet\,\,\circ\,\,\bullet\,\,\bullet\,\,\bullet\,\,\bullet\,\,\bullet\,\,\bullet\,\,\bullet\,\,\bullet\,\,\bullet\,\,\circ\,\,\bullet\,\,\bullet\,\,\bullet\,\,\bullet\,\,\bullet\,\,\bullet\,\,\bullet\,\,\bullet\,\,\bullet\,\,\bullet\,\,\bullet\,\,\bullet\,\,\bullet\,\,\bullet\,\,}
_{\circ\,\,\circ\,\,\circ\,\,\bullet\,\,\circ\,\,\circ\,\,\bullet\,\,\circ\,\,\bullet\,\,\bullet\,\,\circ\,\,\circ\,\,\bullet\,\,\circ\,\,\bullet\,\,\bullet\,\,\circ\,\,\bullet\,\,\bullet\,\,\bullet\,\,\circ\,\,\circ\,\,\bullet\,\,\circ\,\,\bullet\,\,\bullet\,\,\circ\,\,\bullet\,\,\bullet\,\,\bullet\,\,\circ\,\,\bullet\,\,\bullet\,\,\bullet\,\,\bullet\,\,}\\
^{\circ\,\,\circ\,\,\circ\,\,\circ\,\,\bullet\,\,\bullet\,\,\bullet\,\,\bullet\,\,\bullet\,\,\bullet\,\,\bullet\,\,\bullet\,\,\bullet\,\,\bullet\,\,\bullet\,\,\bullet\,\,\bullet\,\,\bullet\,\,\bullet\,\,\bullet\,\,\bullet\,\,\bullet\,\,\bullet\,\,\bullet\,\,\bullet\,\,\bullet\,\,\bullet\,\,\bullet\,\,\bullet\,\,\bullet\,\,\bullet\,\,\bullet\,\,\bullet\,\,\bullet\,\,\bullet\,\,}
_{\circ\,\,\circ\,\,\circ\,\,\circ\,\,\circ\,\,\bullet\,\,\bullet\,\,\bullet\,\,\bullet\,\,\bullet\,\,\circ\,\,\bullet\,\,\bullet\,\,\bullet\,\,\bullet\,\,\bullet\,\,\bullet\,\,\bullet\,\,\bullet\,\,\bullet\,\,\circ\,\,\bullet\,\,\bullet\,\,\bullet\,\,\bullet\,\,\bullet\,\,\bullet\,\,\bullet\,\,\bullet\,\,\bullet\,\,\bullet\,\,\bullet\,\,\bullet\,\,\bullet\,\,\bullet\,\,}\\
^{\circ\,\,\circ\,\,\circ\,\,\circ\,\,\circ\,\,\circ\,\,\bullet\,\,\circ\,\,\bullet\,\,\bullet\,\,\circ\,\,\circ\,\,\bullet\,\,\circ\,\,\bullet\,\,\bullet\,\,\circ\,\,\bullet\,\,\bullet\,\,\bullet\,\,\circ\,\,\circ\,\,\bullet\,\,\circ\,\,\bullet\,\,\bullet\,\,\circ\,\,\bullet\,\,\bullet\,\,\bullet\,\,\circ\,\,\bullet\,\,\bullet\,\,\bullet\,\,\bullet\,\,}
_{\circ\,\,\circ\,\,\circ\,\,\circ\,\,\circ\,\,\circ\,\,\circ\,\,\bullet\,\,\bullet\,\,\bullet\,\,\circ\,\,\circ\,\,\circ\,\,\bullet\,\,\bullet\,\,\bullet\,\,\bullet\,\,\bullet\,\,\bullet\,\,\bullet\,\,\circ\,\,\circ\,\,\circ\,\,\bullet\,\,\bullet\,\,\bullet\,\,\bullet\,\,\bullet\,\,\bullet\,\,\bullet\,\,\bullet\,\,\bullet\,\,\bullet\,\,\bullet\,\,\bullet\,\,}\\
^{\circ\,\,\circ\,\,\circ\,\,\circ\,\,\circ\,\,\circ\,\,\circ\,\,\circ\,\,\bullet\,\,\bullet\,\,\circ\,\,\circ\,\,\circ\,\,\circ\,\,\bullet\,\,\bullet\,\,\circ\,\,\bullet\,\,\bullet\,\,\bullet\,\,\circ\,\,\circ\,\,\circ\,\,\circ\,\,\bullet\,\,\bullet\,\,\circ\,\,\bullet\,\,\bullet\,\,\bullet\,\,\circ\,\,\bullet\,\,\bullet\,\,\bullet\,\,\bullet\,\,}
_{\circ\,\,\circ\,\,\circ\,\,\circ\,\,\circ\,\,\circ\,\,\circ\,\,\circ\,\,\circ\,\,\bullet\,\,\circ\,\,\circ\,\,\circ\,\,\circ\,\,\circ\,\,\bullet\,\,\circ\,\,\circ\,\,\bullet\,\,\bullet\,\,\circ\,\,\circ\,\,\circ\,\,\circ\,\,\circ\,\,\bullet\,\,\circ\,\,\circ\,\,\bullet\,\,\bullet\,\,\circ\,\,\circ\,\,\bullet\,\,\bullet\,\,\bullet\,\,}\\
^{\circ\,\,\circ\,\,\circ\,\,\circ\,\,\circ\,\,\circ\,\,\circ\,\,\circ\,\,\circ\,\,\circ\,\,\bullet\,\,\bullet\,\,\bullet\,\,\bullet\,\,\bullet\,\,\bullet\,\,\bullet\,\,\bullet\,\,\bullet\,\,\bullet\,\,\bullet\,\,\bullet\,\,\bullet\,\,\bullet\,\,\bullet\,\,\bullet\,\,\bullet\,\,\bullet\,\,\bullet\,\,\bullet\,\,\bullet\,\,\bullet\,\,\bullet\,\,\bullet\,\,\bullet\,\,}
_{\circ\,\,\circ\,\,\circ\,\,\circ\,\,\circ\,\,\circ\,\,\circ\,\,\circ\,\,\circ\,\,\circ\,\,\circ\,\,\bullet\,\,\bullet\,\,\bullet\,\,\bullet\,\,\bullet\,\,\bullet\,\,\bullet\,\,\bullet\,\,\bullet\,\,\circ\,\,\bullet\,\,\bullet\,\,\bullet\,\,\bullet\,\,\bullet\,\,\bullet\,\,\bullet\,\,\bullet\,\,\bullet\,\,\bullet\,\,\bullet\,\,\bullet\,\,\bullet\,\,\bullet\,\,}\\
^{\circ\,\,\circ\,\,\circ\,\,\circ\,\,\circ\,\,\circ\,\,\circ\,\,\circ\,\,\circ\,\,\circ\,\,\circ\,\,\circ\,\,\bullet\,\,\circ\,\,\bullet\,\,\bullet\,\,\circ\,\,\bullet\,\,\bullet\,\,\bullet\,\,\circ\,\,\circ\,\,\bullet\,\,\circ\,\,\bullet\,\,\bullet\,\,\circ\,\,\bullet\,\,\bullet\,\,\bullet\,\,\circ\,\,\bullet\,\,\bullet\,\,\bullet\,\,\bullet\,\,}
_{\circ\,\,\circ\,\,\circ\,\,\circ\,\,\circ\,\,\circ\,\,\circ\,\,\circ\,\,\circ\,\,\circ\,\,\circ\,\,\circ\,\,\circ\,\,\bullet\,\,\bullet\,\,\bullet\,\,\bullet\,\,\bullet\,\,\bullet\,\,\bullet\,\,\circ\,\,\circ\,\,\circ\,\,\bullet\,\,\bullet\,\,\bullet\,\,\bullet\,\,\bullet\,\,\bullet\,\,\bullet\,\,\bullet\,\,\bullet\,\,\bullet\,\,\bullet\,\,\bullet\,\,}\\
^{\circ\,\,\circ\,\,\circ\,\,\circ\,\,\circ\,\,\circ\,\,\circ\,\,\circ\,\,\circ\,\,\circ\,\,\circ\,\,\circ\,\,\circ\,\,\circ\,\,\bullet\,\,\bullet\,\,\circ\,\,\bullet\,\,\bullet\,\,\bullet\,\,\circ\,\,\circ\,\,\circ\,\,\circ\,\,\bullet\,\,\bullet\,\,\circ\,\,\bullet\,\,\bullet\,\,\bullet\,\,\circ\,\,\bullet\,\,\bullet\,\,\bullet\,\,\bullet\,\,}
_{\circ\,\,\circ\,\,\circ\,\,\circ\,\,\circ\,\,\circ\,\,\circ\,\,\circ\,\,\circ\,\,\circ\,\,\circ\,\,\circ\,\,\circ\,\,\circ\,\,\circ\,\,\bullet\,\,\circ\,\,\circ\,\,\bullet\,\,\bullet\,\,\circ\,\,\circ\,\,\circ\,\,\circ\,\,\circ\,\,\bullet\,\,\circ\,\,\circ\,\,\bullet\,\,\bullet\,\,\circ\,\,\circ\,\,\bullet\,\,\bullet\,\,\bullet\,\,}\\

^{\circ\,\,\circ\,\,\circ\,\,\circ\,\,\circ\,\,\circ\,\,\circ\,\,\circ\,\,\circ\,\,\circ\,\,\circ\,\,\circ\,\,\circ\,\,\circ\,\,\circ\,\,\circ\,\,\bullet\,\,\bullet\,\,\bullet\,\,\bullet\,\,\circ\,\,\circ\,\,\circ\,\,\circ\,\,\circ\,\,\circ\,\,\bullet\,\,\bullet\,\,\bullet\,\,\bullet\,\,\bullet\,\,\bullet\,\,\bullet\,\,\bullet\,\,\bullet\,\,}
_{\circ\,\,\circ\,\,\circ\,\,\circ\,\,\circ\,\,\circ\,\,\circ\,\,\circ\,\,\circ\,\,\circ\,\,\circ\,\,\circ\,\,\circ\,\,\circ\,\,\circ\,\,\circ\,\,\circ\,\,\bullet\,\,\bullet\,\,\bullet\,\,\circ\,\,\circ\,\,\circ\,\,\circ\,\,\circ\,\,\circ\,\,\circ\,\,\bullet\,\,\bullet\,\,\bullet\,\,\circ\,\,\bullet\,\,\bullet\,\,\bullet\,\,\bullet\,\,}\\

^{\circ\,\,\circ\,\,\circ\,\,\circ\,\,\circ\,\,\circ\,\,\circ\,\,\circ\,\,\circ\,\,\circ\,\,\circ\,\,\circ\,\,\circ\,\,\circ\,\,\circ\,\,\circ\,\,\circ\,\,\circ\,\,\bullet\,\,\bullet\,\,\circ\,\,\circ\,\,\circ\,\,\circ\,\,\circ\,\,\circ\,\,\circ\,\,\circ\,\,\bullet\,\,\bullet\,\,\circ\,\,\circ\,\,\bullet\,\,\bullet\,\,\bullet\,\,}
_{\circ\,\,\circ\,\,\circ\,\,\circ\,\,\circ\,\,\circ\,\,\circ\,\,\circ\,\,\circ\,\,\circ\,\,\circ\,\,\circ\,\,\circ\,\,\circ\,\,\circ\,\,\circ\,\,\circ\,\,\circ\,\,\circ\,\,\bullet\,\,\circ\,\,\circ\,\,\circ\,\,\circ\,\,\circ\,\,\circ\,\,\circ\,\,\circ\,\,\circ\,\,\bullet\,\,\circ\,\,\circ\,\,\circ\,\,\bullet\,\,\bullet\,\,}\\

^{\circ\,\,\circ\,\,\circ\,\,\circ\,\,\circ\,\,\circ\,\,\circ\,\,\circ\,\,\circ\,\,\circ\,\,\circ\,\,\circ\,\,\circ\,\,\circ\,\,\circ\,\,\circ\,\,\circ\,\,\circ\,\,\circ\,\,\circ\,\,\bullet\,\,\bullet\,\,\bullet\,\,\bullet\,\,\bullet\,\,\bullet\,\,\bullet\,\,\bullet\,\,\bullet\,\,\bullet\,\,\bullet\,\,\bullet\,\,\bullet\,\,\bullet\,\,\bullet\,\,}
_{\circ\,\,\circ\,\,\circ\,\,\circ\,\,\circ\,\,\circ\,\,\circ\,\,\circ\,\,\circ\,\,\circ\,\,\circ\,\,\circ\,\,\circ\,\,\circ\,\,\circ\,\,\circ\,\,\circ\,\,\circ\,\,\circ\,\,\circ\,\,\circ\,\,\bullet\,\,\bullet\,\,\bullet\,\,\bullet\,\,\bullet\,\,\bullet\,\,\bullet\,\,\bullet\,\,\bullet\,\,\bullet\,\,\bullet\,\,\bullet\,\,\bullet\,\,\bullet\,\,}\\
^{\circ\,\,\circ\,\,\circ\,\,\circ\,\,\circ\,\,\circ\,\,\circ\,\,\circ\,\,\circ\,\,\circ\,\,\circ\,\,\circ\,\,\circ\,\,\circ\,\,\circ\,\,\circ\,\,\circ\,\,\circ\,\,\circ\,\,\circ\,\,\circ\,\,\circ\,\,\bullet\,\,\circ\,\,\bullet\,\,\bullet\,\,\circ\,\,\bullet\,\,\bullet\,\,\bullet\,\,\circ\,\,\bullet\,\,\bullet\,\,\bullet\,\,\bullet\,\,}
_{\circ\,\,\circ\,\,\circ\,\,\circ\,\,\circ\,\,\circ\,\,\circ\,\,\circ\,\,\circ\,\,\circ\,\,\circ\,\,\circ\,\,\circ\,\,\circ\,\,\circ\,\,\circ\,\,\circ\,\,\circ\,\,\circ\,\,\circ\,\,\circ\,\,\circ\,\,\circ\,\,\bullet\,\,\bullet\,\,\bullet\,\,\bullet\,\,\bullet\,\,\bullet\,\,\bullet\,\,\bullet\,\,\bullet\,\,\bullet\,\,\bullet\,\,\bullet\,\,}\\
^{\circ\,\,\circ\,\,\circ\,\,\circ\,\,\circ\,\,\circ\,\,\circ\,\,\circ\,\,\circ\,\,\circ\,\,\circ\,\,\circ\,\,\circ\,\,\circ\,\,\circ\,\,\circ\,\,\circ\,\,\circ\,\,\circ\,\,\circ\,\,\circ\,\,\circ\,\,\circ\,\,\circ\,\,\bullet\,\,\bullet\,\,\circ\,\,\bullet\,\,\bullet\,\,\bullet\,\,\circ\,\,\bullet\,\,\bullet\,\,\bullet\,\,\bullet\,\,}
_{\circ\,\,\circ\,\,\circ\,\,\circ\,\,\circ\,\,\circ\,\,\circ\,\,\circ\,\,\circ\,\,\circ\,\,\circ\,\,\circ\,\,\circ\,\,\circ\,\,\circ\,\,\circ\,\,\circ\,\,\circ\,\,\circ\,\,\circ\,\,\circ\,\,\circ\,\,\circ\,\,\circ\,\,\circ\,\,\bullet\,\,\circ\,\,\circ\,\,\bullet\,\,\bullet\,\,\circ\,\,\circ\,\,\bullet\,\,\bullet\,\,\bullet\,\,}\\
^{\circ\,\,\circ\,\,\circ\,\,\circ\,\,\circ\,\,\circ\,\,\circ\,\,\circ\,\,\circ\,\,\circ\,\,\circ\,\,\circ\,\,\circ\,\,\circ\,\,\circ\,\,\circ\,\,\circ\,\,\circ\,\,\circ\,\,\circ\,\,\circ\,\,\circ\,\,\circ\,\,\circ\,\,\circ\,\,\circ\,\,\bullet\,\,\bullet\,\,\bullet\,\,\bullet\,\,\bullet\,\,\bullet\,\,\bullet\,\,\bullet\,\,\bullet\,\,}
_{\circ\,\,\circ\,\,\circ\,\,\circ\,\,\circ\,\,\circ\,\,\circ\,\,\circ\,\,\circ\,\,\circ\,\,\circ\,\,\circ\,\,\circ\,\,\circ\,\,\circ\,\,\circ\,\,\circ\,\,\circ\,\,\circ\,\,\circ\,\,\circ\,\,\circ\,\,\circ\,\,\circ\,\,\circ\,\,\circ\,\,\circ\,\,\bullet\,\,\bullet\,\,\bullet\,\,\circ\,\,\bullet\,\,\bullet\,\,\bullet\,\,\bullet\,\,}\\
^{\circ\,\,\circ\,\,\circ\,\,\circ\,\,\circ\,\,\circ\,\,\circ\,\,\circ\,\,\circ\,\,\circ\,\,\circ\,\,\circ\,\,\circ\,\,\circ\,\,\circ\,\,\circ\,\,\circ\,\,\circ\,\,\circ\,\,\circ\,\,\circ\,\,\circ\,\,\circ\,\,\circ\,\,\circ\,\,\circ\,\,\circ\,\,\circ\,\,\bullet\,\,\bullet\,\,\circ\,\,\circ\,\,\bullet\,\,\bullet\,\,\bullet\,\,}
_{\circ\,\,\circ\,\,\circ\,\,\circ\,\,\circ\,\,\circ\,\,\circ\,\,\circ\,\,\circ\,\,\circ\,\,\circ\,\,\circ\,\,\circ\,\,\circ\,\,\circ\,\,\circ\,\,\circ\,\,\circ\,\,\circ\,\,\circ\,\,\circ\,\,\circ\,\,\circ\,\,\circ\,\,\circ\,\,\circ\,\,\circ\,\,\circ\,\,\circ\,\,\bullet\,\,\circ\,\,\circ\,\,\circ\,\,\bullet\,\,\bullet\,\,}\\
^{\circ\,\,\circ\,\,\circ\,\,\circ\,\,\circ\,\,\circ\,\,\circ\,\,\circ\,\,\circ\,\,\circ\,\,\circ\,\,\circ\,\,\circ\,\,\circ\,\,\circ\,\,\circ\,\,\circ\,\,\circ\,\,\circ\,\,\circ\,\,\circ\,\,\circ\,\,\circ\,\,\circ\,\,\circ\,\,\circ\,\,\circ\,\,\circ\,\,\circ\,\,\circ\,\,\bullet\,\,\bullet\,\,\bullet\,\,\bullet\,\,\bullet\,\,}
_{\circ\,\,\circ\,\,\circ\,\,\circ\,\,\circ\,\,\circ\,\,\circ\,\,\circ\,\,\circ\,\,\circ\,\,\circ\,\,\circ\,\,\circ\,\,\circ\,\,\circ\,\,\circ\,\,\circ\,\,\circ\,\,\circ\,\,\circ\,\,\circ\,\,\circ\,\,\circ\,\,\circ\,\,\circ\,\,\circ\,\,\circ\,\,\circ\,\,\circ\,\,\circ\,\,\circ\,\,\bullet\,\,\bullet\,\,\bullet\,\,\bullet\,\,}\\
^{\circ\,\,\circ\,\,\circ\,\,\circ\,\,\circ\,\,\circ\,\,\circ\,\,\circ\,\,\circ\,\,\circ\,\,\circ\,\,\circ\,\,\circ\,\,\circ\,\,\circ\,\,\circ\,\,\circ\,\,\circ\,\,\circ\,\,\circ\,\,\circ\,\,\circ\,\,\circ\,\,\circ\,\,\circ\,\,\circ\,\,\circ\,\,\circ\,\,\circ\,\,\circ\,\,\circ\,\,\circ\,\,\bullet\,\,\bullet\,\,\bullet\,\,}
_{\circ\,\,\circ\,\,\circ\,\,\circ\,\,\circ\,\,\circ\,\,\circ\,\,\circ\,\,\circ\,\,\circ\,\,\circ\,\,\circ\,\,\circ\,\,\circ\,\,\circ\,\,\circ\,\,\circ\,\,\circ\,\,\circ\,\,\circ\,\,\circ\,\,\circ\,\,\circ\,\,\circ\,\,\circ\,\,\circ\,\,\circ\,\,\circ\,\,\circ\,\,\circ\,\,\circ\,\,\circ\,\,\circ\,\,\bullet\,\,\bullet\,\,}\\
^{\circ\,\,\circ\,\,\circ\,\,\circ\,\,\circ\,\,\circ\,\,\circ\,\,\circ\,\,\circ\,\,\circ\,\,\circ\,\,\circ\,\,\circ\,\,\circ\,\,\circ\,\,\circ\,\,\circ\,\,\circ\,\,\circ\,\,\circ\,\,\circ\,\,\circ\,\,\circ\,\,\circ\,\,\circ\,\,\circ\,\,\circ\,\,\circ\,\,\circ\,\,\circ\,\,\circ\,\,\circ\,\,\circ\,\,\circ\,\,\bullet\,\,}
\end{bmatrix}&\\\hline\end{array}$$}

\begin{remark} Dyck words of length $2k$ are bijectively representing the vertex $\mathbb{Z}_{2k+1}$-classes of odd graph $O_k$ \cite{Dejter0}. They are also bijectively representing the vertex $\mathbb{D}_{2k+1}$-classes of middle-levels graph $M_k$ \cite{Dejter2} via lexical matchings \cite{KT}. They allowed a reinterpretation of the Middle-Levels Theorem \cite{gmn,M,MN1,MN2} via nested castling controlled by restricted growth strings (or rgs's, represented in the lower-right of Fig.~3) \cite{Dejter1} yielding a numeral system different from that of item (u) of exercise 6.19 in \cite[p.~222]{Stanley} in that it is a natural topological sorting of  the (ordered) tree of rgs's.
That numeral system also applies to the odd graphs $O_k$ \cite{Dejter0} and their Hamilton cycles via uniform 2-factors formed by $(2k+1)$-cycles \cite{Hcs}.\end{remark}

{\center{\bf TABLE VI}
$$\begin{array}{|l||r|rrrrr|}\hline
k=8&q=1&s_1^1=1;&s_2^1=1,2;&s_3^1=1,2,3;&s_4^1=1,2,3,4;&s_5^1=1,2,3,4,5;\\
i=4&q=2&1,&3,&6,&10,&15;\\
k-i-1=3&q=3&&&&&35.\\\hline
\end{array}$$}

 Let $1<i\in\mathbb{Z}$. Let ${\cal M}_i=\{m_{pq}|1\le p\le i;1\le q\le i\}$ be the $i\times i$-matrix such that: $m_{pq}=1$, if $p\le q$, and $m_{pq}=0$, otherwise.
 For example, $T_2^3={\cal M}_2$, $T_2^4=T_3^4={\cal M}_3$ and $T_2^5=T_4^5={\cal M}_4$. 

 \begin{obs}\label{obs1}
 Let $k>2$. Then $T_1^k=T_k^k=[1]$ and $T_2^k=T_{k-1}^k={\cal M}_{k-1}$.
 \end{obs}

{\center{\bf TABLE VII}
$$\begin{array}{|l||r|rrrrrrr|}\hline
k=11&q=1&s_1^1;&s_1^1,s_2^1;&s_1^1,s_2^1,s_3^1;&s_1^1,...,s_4^1;&s_1^1,...,s_5^1;&s_1^1,...,s_6^1;&s_1^1,...,s_7^1;\\
i=6&q=2&1;&1,   3;&1,   3,      6;&1,...,       10;&1,...,        15;&1,...,       21;&1,...,                28;\\
k-i-1=4&q=3&1,&4,&10,&20,&35,&56,&84;\\
&q=4&&&&&&&210.\\\hline
\end{array}$$}

In view of Observation~\ref{obs1} and the information contained in Table II,  to determine all ordered trees with $k\le 6$ edges, we need the matrices $T_3^6$ and $T_4^6$ of Table IV.
 The entries above the main diagonals of both matrices satisfy the correspondence $t_3^6(p,q)\leftrightarrow t_4^6(11-q,11-p)$. The nonzero entries of $T_3^6$ are precisely
 $t_3^6(3,4)$,
 $t_3^6(3,7)$,
 $t_3^6(5,7)$,
 $t_3^6(6,7)$ and
 $t_3^6(6,8)$, from which the correspondence of entries implied by the second part of Observation~\ref{obs1} yields the respective nonzero entries of $T_4^6$, namely
 $t_4^6(7,8)$,
 $t_4^6(4,8)$,
 $t_4^6(4,6)$,
 $t_4^6(4,5)$ and
 $t_4^6(3,5)$.

\section{Distribution of ordered trees in the tensors $T_i^k$}

{\center{\bf TABLE VIII}
$$\begin{array}{|c|ccccccccccc|ccccccccccc|}\hline
  1&&&&&&&&&&\!\sigma_1^1\!&&&&&&&&&&&\!1\!&\\
  2&&&&&&&&&\!\sigma_1^2\!&&\!\sigma_2^1\!&&&&&&&&&\!1\!&&\!2\!\\
  3&&&&&&&&\!\sigma_1^3\!&&\!\sigma_2^2\!&&&&&&&&&\!1\!&&\!3\!&\\
  4&&&&&&&\!\sigma_1^4\!&&\!\sigma_2^3\!&&\!\sigma_3^2\!&&&&&&&\!1\!&&\!4\!&&\!6\!\\
  5&&&&&&\!\sigma_1^5\!&&\!\sigma_2^4\!&&\!\sigma_3^3\!&&&&&&&\!1\!&&\!5\!&&\!10\!&\\
  6&&&&&\!\sigma_1^6\!&&\!\sigma_2^5\!&&\!\sigma_3^4\!&&\!\sigma_4^3\!&&&&&\!1\!&&\!6\!&&\!15\!&&\!20\!\\
  7&&&&\!\sigma_1^7\!&&\!\sigma_2^6\!&&\!\sigma_3^5\!&&\!\sigma_4^4\!&&&&&\!1\!&&\!7\!&&\!21\!&&\!35\!&\\
  8&&&\!\sigma_1^8\!&&\!\sigma_2^7\!&&\!\sigma_3^6\!&&\!\sigma_4^5\!&&\!\sigma_5^4\!&&&\!1\!&&\!8\!&&\!28\!&&\!56\!&&\!70\!\\
  9&&\!\sigma_1^9\!&&\!\sigma_2^8\!&&\!\sigma_3^7\!&&\!\sigma_4^6\!&&\!\sigma_5^5\!&&&\!1\!&&\!9\!&&\!36\!&&\!84\!&&\!126\!&\\
10&\!\sigma_1^{10}\!&&\!\sigma_2^9\!&&\!\sigma_3^8\!&&\!\sigma_4^7\!&&\!\sigma_5^6\!&&\!\sigma_6^5\!&\!1\!&&\!10\!&&\!45\!&&\!120\!&&\!210\!&&\!252\!\\\hline
\end{array}$$}

For a better visualization, the tensors $T_i^k$ are redrawn
replacing the ones by $\bullet$'s and the zeros by $\circ$'s, as in Table V for $T_3^5$, $T_3^6$, $T_4^6$, $T_3^7$, $T_4^5$, $T_5^7$, $T_3^8$ and $T_4^8$, and
as in Table XI for $T_4^9$.

{\center{\bf TABLE IX}
$$\begin{array}{|c|l|rrrrrrrrr|}\hline
      &q\downarrow|\;p\rightarrow  &1&2   &3	    &4		&5		&6&7		& 8	&9\\\hline
\mbox{A 45}^\circ &1&1&1    &     1    &     1    &     1    &     1    &     1    &     1    &     1  \\
\mbox{rotated}&2&1&2    &     3    &     4    &     5    &     6    &     7    &     8    &     9  \\
\mbox{view}    &3&1&3    &     6    &   10    &   15    &   21    &   28    &   36    &   45  \\
\mbox{of}        &4&1&4    &   10    &   20    &   35    &   56    &   84    & 120    & 165  \\
\mbox{Pascal's}&5&1&5    &   15    &   35    &   70    & 126    & 210    & 330    & 495  \\
\mbox{triangle}&6&1&6    &   21    &   56    & 126    & 252    & 462    & 792    &1167 \\
\Delta&7&1&7    &   28    &   84    & 210    &462     & 924    &1716   &2883 \\
&8&1&8    &   36    & 120    &  330   &792     &1716   &3432   &6315 \\
&9&1&9    &   45    & 165    &  495   &1167   &2883   &6315   &12630\\\hline
\end{array}$$}

\begin{remark}\label{rem} The tensor $T_i^k$ can be presented as a ${k-1\choose i-1}\times{k-1\choose i-1}$-grid ${\cal G}_i^k$ (see Fig.~1 for $(k,i)=(8,4)$), where the $p$-th row header is $A(p)$ and the $q$-th column header is $B(q)^\tau$, for $1\le p,q\le{k-1\choose k-i}$, with $\tau$ meaning ``transpose'', the 0-entries of $T_i^k$ given as light-gray single squares and the 1-entries of $T_i^k$ given as white (i.e. uncolored) single squares.

{\center{\bf TABLE X}
$$\begin{array}{c|cccccccccc}
       &^{5111}&^{3311}&^{3131}&^{2411}&^{2231}&^{2141}&^{1511}&^{1331}&^{1241}&^{1151}\\ \hline
_{5112}&-&2&1&2&1&1&2&1&1&1\\
_{4212}&-&-&1&-&1&1&-&1&1&1\\
_{4113}&-&-&-&5&-&2&5&-&2&2\\
_{3312}&-&-&-&-&1&1&-&1&1&1\\
_{3213}&-&-&-&-&-&2&-&-&2&2\\
_{3114}&-&-&-&-&-&-&9&-&-&3\\
_{2412}&-&-&-&-&-&-&-&1&1&1\\
_{2313}&-&-&-&-&-&-&-&-&2&2\\
_{2214}&-&-&-&-&-&-&-&-&-&3\\
_{2115}&34&-&-&-&-&-&-&-&-&-\\
\end{array}$$}

Apart from the top and rightmost thick black border segments
in ${\cal G}_i^k$, each other vertical (resp. horizontal) segment is: {\bf(a)} thick black if it delimits on the left (resp. bottom) a column (resp. row) headed by an $i$-tuple $B(q)^\tau$ (resp. $A(q)$) with its two last (resp. two penultimate) entries equal to 1; {\bf(b)} thick dark-gray if it delimits on the left (resp. bottom) a column (resp. row) headed by an $i$-tuple $B(q)^\tau$ (resp. $A(q)$) that is not as in (a) but with its last (resp. penultimate) entry equal to 1; {\bf(c)} thin dark-gray, otherwise.\end{remark}

The minimum sub-grids in ${\cal G}_i^k$, as in ${\cal G}_4^8$ in Fig.~1, delimited by thick edges will be called 1{\it-blocks} (resp. 2{\it-blocks}) if those thick edges are dark-gray and/or black (resp. just thick black). Some entries of the headers in Fig.~1 are underlined, namely the first (resp. second) entry below (resp. to the right of) a thick black (resp. dark-gray) line, i.e. at the upper-left corner of a 2-block (resp. 1-block). 

For larger $i$'s, we define $q$-{\it blocks} in ${\cal G}_i^k$ itself, for $q=1,\ldots,i-2$, as follows.
A $q$-block is
a sub-grid of ${\cal G}_i^k$ with the horizontal and vertical headers of its upper-left (resp. lower-right) corner having (resp. just preceding) the first appearance of an integer $\alpha$
 (resp. $\alpha-1$) in the $(k-q-1)$-entry.
 We also say that ${\cal G}_i^k$ is a $(i-1)${\it-block}.

\begin{figure}[htp]
\includegraphics[scale=0.4]{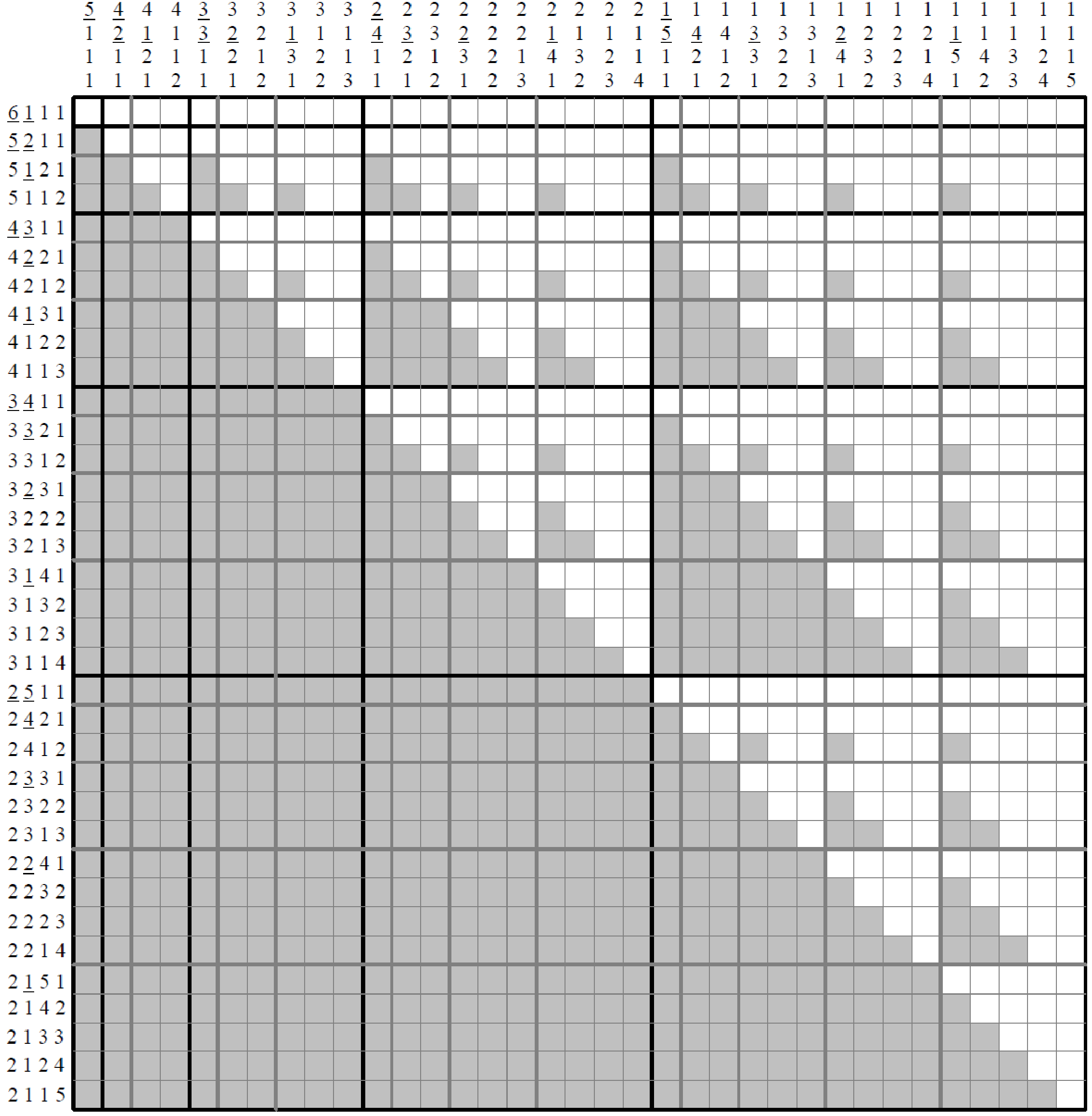}
\caption{${\cal G}_4^8$}
\end{figure}

We will say that the horizontal (resp. vertical) thick lines of ${\cal G}_i^k$, as in Fig.~1, each preceding a header with at least one underlined entry, divide ${\cal G}_i^k$ into horizontal (resp. vertical) $q${\it-strips} if they delimit $q$-blocks. The {\it height} (resp. {\it width}) of a horizontal (resp. vertical) $q$-strip is the number of single squares separating its delimiting thick lines. In ${\cal G}_i^k$, the left-to-right (top-to-bottom) sequences of widths (heights) of vertical (horizontal) $q$-strips are (in both cases) as illustrated in Table VI, where $s_t^1$ stands for the sequence $1,2,\ldots,t$ ($1\le t\in\mathbb{Z}$).

{\center{\bf TABLE XI}
$$\begin{array}{c}
\begin{bmatrix}
^{\bullet\,\,\bullet\,\,\bullet\,\,\bullet\,\,\bullet\,\,\bullet\,\,\bullet\,\,\bullet\,\,\bullet\,\,\bullet\,\,\bullet\,\,\bullet\,\,\bullet\,\,\bullet\,\,\bullet\,\,\bullet\,\,\bullet\,\,\bullet\,\,\bullet\,\,\bullet\,\,\bullet\,\,\bullet\,\,\bullet\,\,\bullet\,\,\bullet\,\,\bullet\,\,\bullet\,\,\bullet\,\,\bullet\,\,\bullet\,\,\bullet\,\,\bullet\,\,\bullet\,\,\bullet\,\,\bullet\,\,\bullet\,\,\bullet\,\,\bullet\,\,\bullet\,\,\bullet\,\,\bullet\,\,\bullet\,\,\bullet\,\,\bullet\,\,\bullet\,\,\bullet\,\,\bullet\,\,\bullet\,\,\bullet\,\,\bullet\,\,\bullet\,\,\bullet\,\,\bullet\,\,\bullet\,\,\bullet\,\,\bullet\,\,}
_{\circ\,\,\bullet\,\,\bullet\,\,\bullet\,\,\bullet\,\,\bullet\,\,\bullet\,\,\bullet\,\,\bullet\,\,\bullet\,\,\bullet\,\,\bullet\,\,\bullet\,\,\bullet\,\,\bullet\,\,\bullet\,\,\bullet\,\,\bullet\,\,\bullet\,\,\bullet\,\,\bullet\,\,\bullet\,\,\bullet\,\,\bullet\,\,\bullet\,\,\bullet\,\,\bullet\,\,\bullet\,\,\bullet\,\,\bullet\,\,\bullet\,\,\bullet\,\,\bullet\,\,\bullet\,\,\bullet\,\,\bullet\,\,\bullet\,\,\bullet\,\,\bullet\,\,\bullet\,\,\bullet\,\,\bullet\,\,\bullet\,\,\bullet\,\,\bullet\,\,\bullet\,\,\bullet\,\,\bullet\,\,\bullet\,\,\bullet\,\,\bullet\,\,\bullet\,\,\bullet\,\,\bullet\,\,\bullet\,\,\bullet\,\,}\\
^{\circ\,\,\circ\,\,\bullet\,\,\bullet\,\,\circ\,\,\bullet\,\,\bullet\,\,\bullet\,\,\bullet\,\,\bullet\,\,\circ\,\,\bullet\,\,\bullet\,\,\bullet\,\,\bullet\,\,\bullet\,\,\bullet\,\,\bullet\,\,\bullet\,\,\bullet\,\,\circ\,\,\bullet\,\,\bullet\,\,\bullet\,\,\bullet\,\,\bullet\,\,\bullet\,\,\bullet\,\,\bullet\,\,\bullet\,\,\bullet\,\,\bullet\,\,\bullet\,\,\bullet\,\,\bullet\,\,\circ\,\,\bullet\,\,\bullet\,\,\bullet\,\,\bullet\,\,\bullet\,\,\bullet\,\,\bullet\,\,\bullet\,\,\bullet\,\,\bullet\,\,\bullet\,\,\bullet\,\,\bullet\,\,\bullet\,\,\bullet\,\,\bullet\,\,\bullet\,\,\bullet\,\,\bullet\,\,\bullet\,\,}
_{\circ\,\,\circ\,\,\circ\,\,\bullet\,\,\circ\,\,\circ\,\,\bullet\,\,\circ\,\,\bullet\,\,\bullet\,\,\circ\,\,\circ\,\,\bullet\,\,\circ\,\,\bullet\,\,\bullet\,\,\circ\,\,\bullet\,\,\bullet\,\,\bullet\,\,\circ\,\,\circ\,\,\bullet\,\,\circ\,\,\bullet\,\,\bullet\,\,\circ\,\,\bullet\,\,\bullet\,\,\bullet\,\,\circ\,\,\bullet\,\,\bullet\,\,\bullet\,\,\bullet\,\,\circ\,\,\circ\,\,\bullet\,\,\circ\,\,\bullet\,\,\bullet\,\,\circ\,\,\bullet\,\,\bullet\,\,\bullet\,\,\circ\,\,\bullet\,\,\bullet\,\,\bullet\,\,\bullet\,\,\circ\,\,\bullet\,\,\bullet\,\,\bullet\,\,\bullet\,\,\bullet\,\,}\\
^{\circ\,\,\circ\,\,\circ\,\,\circ\,\,\bullet\,\,\bullet\,\,\bullet\,\,\bullet\,\,\bullet\,\,\bullet\,\,\bullet\,\,\bullet\,\,\bullet\,\,\bullet\,\,\bullet\,\,\bullet\,\,\bullet\,\,\bullet\,\,\bullet\,\,\bullet\,\,\bullet\,\,\bullet\,\,\bullet\,\,\bullet\,\,\bullet\,\,\bullet\,\,\bullet\,\,\bullet\,\,\bullet\,\,\bullet\,\,\bullet\,\,\bullet\,\,\bullet\,\,\bullet\,\,\bullet\,\,\bullet\,\,\bullet\,\,\bullet\,\,\bullet\,\,\bullet\,\,\bullet\,\,\bullet\,\,\bullet\,\,\bullet\,\,\bullet\,\,\bullet\,\,\bullet\,\,\bullet\,\,\bullet\,\,\bullet\,\,\bullet\,\,\bullet\,\,\bullet\,\,\bullet\,\,\bullet\,\,\bullet\,\,}
_{\circ\,\,\circ\,\,\circ\,\,\circ\,\,\circ\,\,\bullet\,\,\bullet\,\,\bullet\,\,\bullet\,\,\bullet\,\,\circ\,\,\bullet\,\,\bullet\,\,\bullet\,\,\bullet\,\,\bullet\,\,\bullet\,\,\bullet\,\,\bullet\,\,\bullet\,\,\circ\,\,\bullet\,\,\bullet\,\,\bullet\,\,\bullet\,\,\bullet\,\,\bullet\,\,\bullet\,\,\bullet\,\,\bullet\,\,\bullet\,\,\bullet\,\,\bullet\,\,\bullet\,\,\bullet\,\,\circ\,\,\bullet\,\,\bullet\,\,\bullet\,\,\bullet\,\,\bullet\,\,\bullet\,\,\bullet\,\,\bullet\,\,\bullet\,\,\bullet\,\,\bullet\,\,\bullet\,\,\bullet\,\,\bullet\,\,\bullet\,\,\bullet\,\,\bullet\,\,\bullet\,\,\bullet\,\,\bullet\,\,}\\
^{\circ\,\,\circ\,\,\circ\,\,\circ\,\,\circ\,\,\circ\,\,\bullet\,\,\circ\,\,\bullet\,\,\bullet\,\,\circ\,\,\circ\,\,\bullet\,\,\circ\,\,\bullet\,\,\bullet\,\,\circ\,\,\bullet\,\,\bullet\,\,\bullet\,\,\circ\,\,\circ\,\,\bullet\,\,\circ\,\,\bullet\,\,\bullet\,\,\circ\,\,\bullet\,\,\bullet\,\,\bullet\,\,\circ\,\,\bullet\,\,\bullet\,\,\bullet\,\,\bullet\,\,\circ\,\,\circ\,\,\bullet\,\,\circ\,\,\bullet\,\,\bullet\,\,\circ\,\,\bullet\,\,\bullet\,\,\bullet\,\,\circ\,\,\bullet\,\,\bullet\,\,\bullet\,\,\bullet\,\,\circ\,\,\bullet\,\,\bullet\,\,\bullet\,\,\bullet\,\,\bullet\,\,}
_{\circ\,\,\circ\,\,\circ\,\,\circ\,\,\circ\,\,\circ\,\,\circ\,\,\bullet\,\,\bullet\,\,\bullet\,\,\circ\,\,\circ\,\,\circ\,\,\bullet\,\,\bullet\,\,\bullet\,\,\bullet\,\,\bullet\,\,\bullet\,\,\bullet\,\,\circ\,\,\circ\,\,\circ\,\,\bullet\,\,\bullet\,\,\bullet\,\,\bullet\,\,\bullet\,\,\bullet\,\,\bullet\,\,\bullet\,\,\bullet\,\,\bullet\,\,\bullet\,\,\bullet\,\,\circ\,\,\circ\,\,\circ\,\,\bullet\,\,\bullet\,\,\bullet\,\,\bullet\,\,\bullet\,\,\bullet\,\,\bullet\,\,\bullet\,\,\bullet\,\,\bullet\,\,\bullet\,\,\bullet\,\,\bullet\,\,\bullet\,\,\bullet\,\,\bullet\,\,\bullet\,\,\bullet\,\,}\\
^{\circ\,\,\circ\,\,\circ\,\,\circ\,\,\circ\,\,\circ\,\,\circ\,\,\circ\,\,\bullet\,\,\bullet\,\,\circ\,\,\circ\,\,\circ\,\,\circ\,\,\bullet\,\,\bullet\,\,\circ\,\,\bullet\,\,\bullet\,\,\bullet\,\,\circ\,\,\circ\,\,\circ\,\,\circ\,\,\bullet\,\,\bullet\,\,\circ\,\,\bullet\,\,\bullet\,\,\bullet\,\,\circ\,\,\bullet\,\,\bullet\,\,\bullet\,\,\bullet\,\,\circ\,\,\circ\,\,\circ\,\,\circ\,\,\bullet\,\,\bullet\,\,\circ\,\,\bullet\,\,\bullet\,\,\bullet\,\,\circ\,\,\bullet\,\,\bullet\,\,\bullet\,\,\bullet\,\,\circ\,\,\bullet\,\,\bullet\,\,\bullet\,\,\bullet\,\,\bullet\,\,}
_{\circ\,\,\circ\,\,\circ\,\,\circ\,\,\circ\,\,\circ\,\,\circ\,\,\circ\,\,\circ\,\,\bullet\,\,\circ\,\,\circ\,\,\circ\,\,\circ\,\,\circ\,\,\bullet\,\,\circ\,\,\circ\,\,\bullet\,\,\bullet\,\,\circ\,\,\circ\,\,\circ\,\,\circ\,\,\circ\,\,\bullet\,\,\circ\,\,\circ\,\,\bullet\,\,\bullet\,\,\circ\,\,\circ\,\,\bullet\,\,\bullet\,\,\bullet\,\,\circ\,\,\circ\,\,\circ\,\,\circ\,\,\circ\,\,\bullet\,\,\circ\,\,\circ\,\,\bullet\,\,\bullet\,\,\circ\,\,\circ\,\,\bullet\,\,\bullet\,\,\bullet\,\,\circ\,\,\circ\,\,\bullet\,\,\bullet\,\,\bullet\,\,\bullet\,\,}\\
^{\circ\,\,\circ\,\,\circ\,\,\circ\,\,\circ\,\,\circ\,\,\circ\,\,\circ\,\,\circ\,\,\circ\,\,\bullet\,\,\bullet\,\,\bullet\,\,\bullet\,\,\bullet\,\,\bullet\,\,\bullet\,\,\bullet\,\,\bullet\,\,\bullet\,\,\bullet\,\,\bullet\,\,\bullet\,\,\bullet\,\,\bullet\,\,\bullet\,\,\bullet\,\,\bullet\,\,\bullet\,\,\bullet\,\,\bullet\,\,\bullet\,\,\bullet\,\,\bullet\,\,\bullet\,\,\bullet\,\,\bullet\,\,\bullet\,\,\bullet\,\,\bullet\,\,\bullet\,\,\bullet\,\,\bullet\,\,\bullet\,\,\bullet\,\,\bullet\,\,\bullet\,\,\bullet\,\,\bullet\,\,\bullet\,\,\bullet\,\,\bullet\,\,\bullet\,\,\bullet\,\,\bullet\,\,\bullet\,\,}
_{\circ\,\,\circ\,\,\circ\,\,\circ\,\,\circ\,\,\circ\,\,\circ\,\,\circ\,\,\circ\,\,\circ\,\,\circ\,\,\bullet\,\,\bullet\,\,\bullet\,\,\bullet\,\,\bullet\,\,\bullet\,\,\bullet\,\,\bullet\,\,\bullet\,\,\circ\,\,\bullet\,\,\bullet\,\,\bullet\,\,\bullet\,\,\bullet\,\,\bullet\,\,\bullet\,\,\bullet\,\,\bullet\,\,\bullet\,\,\bullet\,\,\bullet\,\,\bullet\,\,\bullet\,\,\circ\,\,\bullet\,\,\bullet\,\,\bullet\,\,\bullet\,\,\bullet\,\,\bullet\,\,\bullet\,\,\bullet\,\,\bullet\,\,\bullet\,\,\bullet\,\,\bullet\,\,\bullet\,\,\bullet\,\,\bullet\,\,\bullet\,\,\bullet\,\,\bullet\,\,\bullet\,\,\bullet\,\,}\\
^{\circ\,\,\circ\,\,\circ\,\,\circ\,\,\circ\,\,\circ\,\,\circ\,\,\circ\,\,\circ\,\,\circ\,\,\circ\,\,\circ\,\,\bullet\,\,\circ\,\,\bullet\,\,\bullet\,\,\circ\,\,\bullet\,\,\bullet\,\,\bullet\,\,\circ\,\,\circ\,\,\bullet\,\,\circ\,\,\bullet\,\,\bullet\,\,\circ\,\,\bullet\,\,\bullet\,\,\bullet\,\,\circ\,\,\bullet\,\,\bullet\,\,\bullet\,\,\bullet\,\,\circ\,\,\circ\,\,\bullet\,\,\circ\,\,\bullet\,\,\bullet\,\,\circ\,\,\bullet\,\,\bullet\,\,\bullet\,\,\circ\,\,\bullet\,\,\bullet\,\,\bullet\,\,\bullet\,\,\circ\,\,\bullet\,\,\bullet\,\,\bullet\,\,\bullet\,\,\bullet\,\,}
_{\circ\,\,\circ\,\,\circ\,\,\circ\,\,\circ\,\,\circ\,\,\circ\,\,\circ\,\,\circ\,\,\circ\,\,\circ\,\,\circ\,\,\circ\,\,\bullet\,\,\bullet\,\,\bullet\,\,\bullet\,\,\bullet\,\,\bullet\,\,\bullet\,\,\circ\,\,\circ\,\,\circ\,\,\bullet\,\,\bullet\,\,\bullet\,\,\bullet\,\,\bullet\,\,\bullet\,\,\bullet\,\,\bullet\,\,\bullet\,\,\bullet\,\,\bullet\,\,\bullet\,\,\circ\,\,\circ\,\,\circ\,\,\bullet\,\,\bullet\,\,\bullet\,\,\bullet\,\,\bullet\,\,\bullet\,\,\bullet\,\,\bullet\,\,\bullet\,\,\bullet\,\,\bullet\,\,\bullet\,\,\bullet\,\,\bullet\,\,\bullet\,\,\bullet\,\,\bullet\,\,\bullet\,\,}\\
^{\circ\,\,\circ\,\,\circ\,\,\circ\,\,\circ\,\,\circ\,\,\circ\,\,\circ\,\,\circ\,\,\circ\,\,\circ\,\,\circ\,\,\circ\,\,\circ\,\,\bullet\,\,\bullet\,\,\circ\,\,\bullet\,\,\bullet\,\,\bullet\,\,\circ\,\,\circ\,\,\circ\,\,\circ\,\,\bullet\,\,\bullet\,\,\circ\,\,\bullet\,\,\bullet\,\,\bullet\,\,\circ\,\,\bullet\,\,\bullet\,\,\bullet\,\,\bullet\,\,\circ\,\,\circ\,\,\circ\,\,\circ\,\,\bullet\,\,\bullet\,\,\circ\,\,\bullet\,\,\bullet\,\,\bullet\,\,\circ\,\,\bullet\,\,\bullet\,\,\bullet\,\,\bullet\,\,\circ\,\,\bullet\,\,\bullet\,\,\bullet\,\,\bullet\,\,\bullet\,\,}
_{\circ\,\,\circ\,\,\circ\,\,\circ\,\,\circ\,\,\circ\,\,\circ\,\,\circ\,\,\circ\,\,\circ\,\,\circ\,\,\circ\,\,\circ\,\,\circ\,\,\circ\,\,\bullet\,\,\circ\,\,\circ\,\,\bullet\,\,\bullet\,\,\circ\,\,\circ\,\,\circ\,\,\circ\,\,\circ\,\,\bullet\,\,\circ\,\,\circ\,\,\bullet\,\,\bullet\,\,\circ\,\,\circ\,\,\bullet\,\,\bullet\,\,\bullet\,\,\circ\,\,\circ\,\,\circ\,\,\circ\,\,\circ\,\,\bullet\,\,\circ\,\,\circ\,\,\bullet\,\,\bullet\,\,\circ\,\,\circ\,\,\bullet\,\,\bullet\,\,\bullet\,\,\circ\,\,\circ\,\,\bullet\,\,\bullet\,\,\bullet\,\,\bullet\,\,}\\
^{\circ\,\,\circ\,\,\circ\,\,\circ\,\,\circ\,\,\circ\,\,\circ\,\,\circ\,\,\circ\,\,\circ\,\,\circ\,\,\circ\,\,\circ\,\,\circ\,\,\circ\,\,\circ\,\,\bullet\,\,\bullet\,\,\bullet\,\,\bullet\,\,\circ\,\,\circ\,\,\circ\,\,\circ\,\,\circ\,\,\circ  \,\,\bullet\,\,\bullet\,\,\bullet\,\,\bullet\,\,\bullet\,\,\bullet\,\,\bullet\,\,\bullet\,\,\bullet\,\,\circ\,\,\circ\,\,\circ\,\,\circ\,\,\circ\,\,\circ\,\,\bullet\,\,\bullet\,\,\bullet\,\,\bullet\,\,\bullet\,\,\bullet\,\,\bullet\,\,\bullet\,\,\bullet\,\,\bullet\,\,\bullet\,\,\bullet\,\,\bullet\,\,\bullet\,\,\bullet\,\,}
_{\circ\,\,\circ\,\,\circ\,\,\circ\,\,\circ\,\,\circ\,\,\circ\,\,\circ\,\,\circ\,\,\circ\,\,\circ\,\,\circ\,\,\circ\,\,\circ\,\,\circ\,\,\circ\,\,\circ\,\,\bullet\,\,\bullet\,\,\bullet\,\,\circ\,\,\circ\,\,\circ\,\,\circ\,\,\circ\,\,\circ  \,\,\circ\,\,\bullet\,\,\bullet\,\,\bullet\,\,\circ\,\,\bullet\,\,\bullet\,\,\bullet\,\,\bullet\,\,\circ\,\,\circ\,\,\circ\,\,\circ\,\,\circ\,\,\circ\,\,\circ\,\,\bullet\,\,\bullet\,\,\bullet\,\,\circ\,\,\bullet\,\,\bullet\,\,\bullet\,\,\bullet\,\,\circ\,\,\bullet\,\,\bullet\,\,\bullet\,\,\bullet\,\,\bullet\,\,}\\
^{\circ\,\,\circ\,\,\circ\,\,\circ\,\,\circ\,\,\circ\,\,\circ\,\,\circ\,\,\circ\,\,\circ\,\,\circ\,\,\circ\,\,\circ\,\,\circ\,\,\circ\,\,\circ\,\,\circ\,\,\circ\,\,\bullet\,\,\bullet\,\,\circ\,\,\circ\,\,\circ\,\,\circ\,\,\circ\,\,\circ  \,\,\circ\,\,\circ\,\,\bullet\,\,\bullet\,\,\circ\,\,\circ\,\,\bullet\,\,\bullet\,\,\bullet\,\,\circ\,\,\circ\,\,\circ\,\,\circ\,\,\circ\,\,\circ\,\,\circ\,\,\circ\,\,\bullet\,\,\bullet\,\,\circ\,\,\circ\,\,\bullet\,\,\bullet\,\,\bullet\,\,\circ\,\,\circ\,\,\bullet\,\,\bullet\,\,\bullet\,\,\bullet\,\,}
_{\circ\,\,\circ\,\,\circ\,\,\circ\,\,\circ\,\,\circ\,\,\circ\,\,\circ\,\,\circ\,\,\circ\,\,\circ\,\,\circ\,\,\circ\,\,\circ\,\,\circ\,\,\circ\,\,\circ\,\,\circ\,\,\circ\,\,\bullet\,\,\circ\,\,\circ\,\,\circ\,\,\circ\,\,\circ\,\,\circ  \,\,\circ\,\,\circ\,\,\circ\,\,\bullet\,\,\circ\,\,\circ\,\,\circ\,\,\bullet\,\,\bullet\,\,\circ\,\,\circ\,\,\circ\,\,\circ\,\,\circ\,\,\circ\,\,\circ\,\,\circ\,\,\circ\,\,\bullet\,\,\circ\,\,\circ\,\,\circ\,\,\bullet\,\,\bullet\,\,\circ\,\,\circ\,\,\circ\,\,\bullet\,\,\bullet\,\,\bullet\,\,}\\
^{\circ\,\,\circ\,\,\circ\,\,\circ\,\,\circ\,\,\circ\,\,\circ\,\,\circ\,\,\circ\,\,\circ\,\,\circ\,\,\circ\,\,\circ\,\,\circ\,\,\circ\,\,\circ\,\,\circ\,\,\circ\,\,\circ\,\,\circ\,\,\bullet\,\,\bullet\,\,\bullet\,\,\bullet\,\,\bullet\,\,\bullet\,\,\bullet\,\,\bullet\,\,\bullet\,\,\bullet\,\,\bullet\,\,\bullet\,\,\bullet\,\,\bullet\,\,\bullet\,\,\bullet\,\,\bullet\,\,\bullet\,\,\bullet\,\,\bullet\,\,\bullet\,\,\bullet\,\,\bullet\,\,\bullet\,\,\bullet\,\,\bullet\,\,\bullet\,\,\bullet\,\,\bullet\,\,\bullet\,\,\bullet\,\,\bullet\,\,\bullet\,\,\bullet\,\,\bullet\,\,\bullet\,\,}
_{\circ\,\,\circ\,\,\circ\,\,\circ\,\,\circ\,\,\circ\,\,\circ\,\,\circ\,\,\circ\,\,\circ\,\,\circ\,\,\circ\,\,\circ\,\,\circ\,\,\circ\,\,\circ\,\,\circ\,\,\circ\,\,\circ\,\,\circ\,\,\circ\,\,\bullet\,\,\bullet\,\,\bullet\,\,\bullet\,\,\bullet\,\,\bullet\,\,\bullet\,\,\bullet\,\,\bullet\,\,\bullet\,\,\bullet\,\,\bullet\,\,\bullet\,\,\bullet\,\,\circ\,\,\bullet\,\,\bullet\,\,\bullet\,\,\bullet\,\,\bullet\,\,\bullet\,\,\bullet\,\,\bullet\,\,\bullet\,\,\bullet\,\,\bullet\,\,\bullet\,\,\bullet\,\,\bullet\,\,\bullet\,\,\bullet\,\,\bullet\,\,\bullet\,\,\bullet\,\,\bullet\,\,}\\
^{\circ\,\,\circ\,\,\circ\,\,\circ\,\,\circ\,\,\circ\,\,\circ\,\,\circ\,\,\circ\,\,\circ\,\,\circ\,\,\circ\,\,\circ\,\,\circ\,\,\circ\,\,\circ\,\,\circ\,\,\circ\,\,\circ\,\,\circ\,\,\circ\,\,\circ\,\,\bullet\,\,\circ\,\,\bullet\,\,\bullet\,\,\circ\,\,\bullet\,\,\bullet\,\,\bullet\,\,\circ\,\,\bullet\,\,\bullet\,\,\bullet\,\,\bullet\,\,\circ\,\,\circ\,\,\bullet\,\,\circ\,\,\bullet\,\,\bullet\,\,\circ\,\,\bullet\,\,\bullet\,\,\bullet\,\,\circ\,\,\bullet\,\,\bullet\,\,\bullet\,\,\bullet\,\,\circ\,\,\bullet\,\,\bullet\,\,\bullet\,\,\bullet\,\,\bullet\,\,}
_{\circ\,\,\circ\,\,\circ\,\,\circ\,\,\circ\,\,\circ\,\,\circ\,\,\circ\,\,\circ\,\,\circ\,\,\circ\,\,\circ\,\,\circ\,\,\circ\,\,\circ\,\,\circ\,\,\circ\,\,\circ\,\,\circ\,\,\circ\,\,\circ\,\,\circ\,\,\circ\,\,\bullet\,\,\bullet\,\,\bullet\,\,\bullet\,\,\bullet\,\,\bullet\,\,\bullet\,\,\bullet\,\,\bullet\,\,\bullet\,\,\bullet\,\,\bullet\,\,\circ\,\,\circ\,\,\circ\,\,\bullet\,\,\bullet\,\,\bullet\,\,\bullet\,\,\bullet\,\,\bullet\,\,\bullet\,\,\bullet\,\,\bullet\,\,\bullet\,\,\bullet\,\,\bullet\,\,\bullet\,\,\bullet\,\,\bullet\,\,\bullet\,\,\bullet\,\,\bullet\,\,}\\
^{\circ\,\,\circ\,\,\circ\,\,\circ\,\,\circ\,\,\circ\,\,\circ\,\,\circ\,\,\circ\,\,\circ\,\,\circ\,\,\circ\,\,\circ\,\,\circ\,\,\circ\,\,\circ\,\,\circ\,\,\circ\,\,\circ\,\,\circ\,\,\circ\,\,\circ\,\,\circ\,\,\circ\,\,\bullet\,\,\bullet\,\,\circ\,\,\bullet\,\,\bullet\,\,\bullet\,\,\circ\,\,\bullet\,\,\bullet\,\,\bullet\,\,\bullet\,\,\circ\,\,\circ\,\,\circ\,\,\circ\,\,\bullet\,\,\bullet\,\,\circ\,\,\bullet\,\,\bullet\,\,\bullet\,\,\circ\,\,\bullet\,\,\bullet\,\,\bullet\,\,\bullet\,\,\circ\,\,\bullet\,\,\bullet\,\,\bullet\,\,\bullet\,\,\bullet\,\,}
_{\circ\,\,\circ\,\,\circ\,\,\circ\,\,\circ\,\,\circ\,\,\circ\,\,\circ\,\,\circ\,\,\circ\,\,\circ\,\,\circ\,\,\circ\,\,\circ\,\,\circ\,\,\circ\,\,\circ\,\,\circ\,\,\circ\,\,\circ\,\,\circ\,\,\circ\,\,\circ\,\,\circ\,\,\circ\,\,\bullet\,\,\circ\,\,\circ\,\,\bullet\,\,\bullet\,\,\circ\,\,\circ\,\,\bullet\,\,\bullet\,\,\bullet\,\,\circ\,\,\circ\,\,\circ\,\,\circ\,\,\circ\,\,\bullet\,\,\circ\,\,\circ\,\,\bullet\,\,\bullet\,\,\circ\,\,\circ\,\,\bullet\,\,\bullet\,\,\bullet\,\,\circ\,\,\circ\,\,\bullet\,\,\bullet\,\,\bullet\,\,\bullet\,\,}\\
^{\circ\,\,\circ\,\,\circ\,\,\circ\,\,\circ\,\,\circ\,\,\circ\,\,\circ\,\,\circ\,\,\circ\,\,\circ\,\,\circ\,\,\circ\,\,\circ\,\,\circ\,\,\circ\,\,\circ\,\,\circ\,\,\circ\,\,\circ\,\,\circ\,\,\circ\,\,\circ\,\,\circ\,\,\circ\,\,\circ\,\,\bullet\,\,\bullet\,\,\bullet\,\,\bullet\,\,\bullet\,\,\bullet\,\,\bullet\,\,\bullet\,\,\bullet\,\,\circ\,\,\circ\,\,\circ\,\,\circ\,\,\circ\,\,\circ\,\,\bullet\,\,\bullet\,\,\bullet\,\,\bullet\,\,\bullet\,\,\bullet\,\,\bullet\,\,\bullet\,\,\bullet\,\,\bullet\,\,\bullet\,\,\bullet\,\,\bullet\,\,\bullet\,\,\bullet\,\,}
_{\circ\,\,\circ\,\,\circ\,\,\circ\,\,\circ\,\,\circ\,\,\circ\,\,\circ\,\,\circ\,\,\circ\,\,\circ\,\,\circ\,\,\circ\,\,\circ\,\,\circ\,\,\circ\,\,\circ\,\,\circ\,\,\circ\,\,\circ\,\,\circ\,\,\circ\,\,\circ\,\,\circ\,\,\circ\,\,\circ\,\,\circ\,\,\bullet\,\,\bullet\,\,\bullet\,\,\circ\,\,\bullet\,\,\bullet\,\,\bullet\,\,\bullet\,\,\circ\,\,\circ\,\,\circ\,\,\circ\,\,\circ\,\,\circ\,\,\circ\,\,\bullet\,\,\bullet\,\,\bullet\,\,\circ\,\,\bullet\,\,\bullet\,\,\bullet\,\,\bullet\,\,\circ\,\,\bullet\,\,\bullet\,\,\bullet\,\,\bullet\,\,\bullet\,\,}\\
^{\circ\,\,\circ\,\,\circ\,\,\circ\,\,\circ\,\,\circ\,\,\circ\,\,\circ\,\,\circ\,\,\circ\,\,\circ\,\,\circ\,\,\circ\,\,\circ\,\,\circ\,\,\circ\,\,\circ\,\,\circ\,\,\circ\,\,\circ\,\,\circ\,\,\circ\,\,\circ\,\,\circ\,\,\circ\,\,\circ\,\,\circ\,\,\circ\,\,\bullet\,\,\bullet\,\,\circ\,\,\circ\,\,\bullet\,\,\bullet\,\,\bullet\,\,\circ\,\,\circ\,\,\circ\,\,\circ\,\,\circ\,\,\circ\,\,\circ\,\,\circ\,\,\bullet\,\,\bullet\,\,\circ\,\,\circ\,\,\bullet\,\,\bullet\,\,\bullet\,\,\circ\,\,\circ\,\,\bullet\,\,\bullet\,\,\bullet\,\,\bullet\,\,}
_{\circ\,\,\circ\,\,\circ\,\,\circ\,\,\circ\,\,\circ\,\,\circ\,\,\circ\,\,\circ\,\,\circ\,\,\circ\,\,\circ\,\,\circ\,\,\circ\,\,\circ\,\,\circ\,\,\circ\,\,\circ\,\,\circ\,\,\circ\,\,\circ\,\,\circ\,\,\circ\,\,\circ\,\,\circ\,\,\circ\,\,\circ\,\,\circ\,\,\circ\,\,\bullet\,\,\circ\,\,\circ\,\,\circ\,\,\bullet\,\,\bullet\,\,\circ\,\,\circ\,\,\circ\,\,\circ\,\,\circ\,\,\circ\,\,\circ\,\,\circ\,\,\circ\,\,\bullet\,\,\circ\,\,\circ\,\,\circ\,\,\bullet\,\,\bullet\,\,\circ\,\,\circ\,\,\circ\,\,\bullet\,\,\bullet\,\,\bullet\,\,}\\
^{\circ\,\,\circ\,\,\circ\,\,\circ\,\,\circ\,\,\circ\,\,\circ\,\,\circ\,\,\circ\,\,\circ\,\,\circ\,\,\circ\,\,\circ\,\,\circ\,\,\circ\,\,\circ\,\,\circ\,\,\circ\,\,\circ\,\,\circ\,\,\circ\,\,\circ\,\,\circ\,\,\circ\,\,\circ\,\,\circ\,\,\circ\,\,\circ\,\,\circ\,\,\circ\,\,\bullet\,\,\bullet\,\,\bullet\,\,\bullet\,\,\bullet\,\,\circ\,\,\circ\,\,\circ\,\,\circ\,\,\circ\,\,\circ\,\,\circ\,\,\circ\,\,\circ\,\,\circ\,\,\bullet\,\,\bullet\,\,\bullet\,\,\bullet\,\,\bullet\,\,\bullet\,\,\bullet\,\,\bullet\,\,\bullet\,\,\bullet\,\,\bullet\,\,}
_{\circ\,\,\circ\,\,\circ\,\,\circ\,\,\circ\,\,\circ\,\,\circ\,\,\circ\,\,\circ\,\,\circ\,\,\circ\,\,\circ\,\,\circ\,\,\circ\,\,\circ\,\,\circ\,\,\circ\,\,\circ\,\,\circ\,\,\circ\,\,\circ\,\,\circ\,\,\circ\,\,\circ\,\,\circ\,\,\circ\,\,\circ\,\,\circ\,\,\circ\,\,\circ\,\,\circ\,\,\bullet\,\,\bullet\,\,\bullet\,\,\bullet\,\,\circ\,\,\circ\,\,\circ\,\,\circ\,\,\circ\,\,\circ\,\,\circ\,\,\circ\,\,\circ\,\,\circ\,\,\circ\,\,\bullet\,\,\bullet\,\,\bullet\,\,\bullet\,\,\circ\,\,\bullet\,\,\bullet\,\,\bullet\,\,\bullet\,\,\bullet\,\,}\\
^{\circ\,\,\circ\,\,\circ\,\,\circ\,\,\circ\,\,\circ\,\,\circ\,\,\circ\,\,\circ\,\,\circ\,\,\circ\,\,\circ\,\,\circ\,\,\circ\,\,\circ\,\,\circ\,\,\circ\,\,\circ\,\,\circ\,\,\circ\,\,\circ\,\,\circ\,\,\circ\,\,\circ\,\,\circ\,\,\circ\,\,\circ\,\,\circ\,\,\circ\,\,\circ\,\,\circ\,\,\circ\,\,\bullet\,\,\bullet\,\,\bullet\,\,\circ\,\,\circ\,\,\circ\,\,\circ\,\,\circ\,\,\circ\,\,\circ\,\,\circ\,\,\circ\,\,\circ\,\,\circ\,\,\circ\,\,\bullet\,\,\bullet\,\,\bullet\,\,\circ\,\,\circ\,\,\bullet\,\,\bullet\,\,\bullet\,\,\bullet\,\,}
_{\circ\,\,\circ\,\,\circ\,\,\circ\,\,\circ\,\,\circ\,\,\circ\,\,\circ\,\,\circ\,\,\circ\,\,\circ\,\,\circ\,\,\circ\,\,\circ\,\,\circ\,\,\circ\,\,\circ\,\,\circ\,\,\circ\,\,\circ\,\,\circ\,\,\circ\,\,\circ\,\,\circ\,\,\circ\,\,\circ\,\,\circ\,\,\circ\,\,\circ\,\,\circ\,\,\circ\,\,\circ\,\,\circ\,\,\bullet\,\,\bullet\,\,\circ\,\,\circ\,\,\circ\,\,\circ\,\,\circ\,\,\circ\,\,\circ\,\,\circ\,\,\circ\,\,\circ\,\,\circ\,\,\circ\,\,\circ\,\,\bullet\,\,\bullet\,\,\circ\,\,\circ\,\,\circ\,\,\bullet\,\,\bullet\,\,\bullet\,\,}\\
^{\circ\,\,\circ\,\,\circ\,\,\circ\,\,\circ\,\,\circ\,\,\circ\,\,\circ\,\,\circ\,\,\circ\,\,\circ\,\,\circ\,\,\circ\,\,\circ\,\,\circ\,\,\circ\,\,\circ\,\,\circ\,\,\circ\,\,\circ\,\,\circ\,\,\circ\,\,\circ\,\,\circ\,\,\circ\,\,\circ\,\,\circ\,\,\circ\,\,\circ\,\,\circ\,\,\circ\,\,\circ\,\,\circ\,\,\circ\,\,\bullet\,\,\circ\,\,\circ\,\,\circ\,\,\circ\,\,\circ\,\,\circ\,\,\circ\,\,\circ\,\,\circ\,\,\circ\,\,\circ\,\,\circ\,\,\circ\,\,\circ\,\,\bullet\,\,\circ\,\,\circ\,\,\circ\,\,\circ\,\,\bullet\,\,\bullet\,\,}
_{\circ\,\,\circ\,\,\circ\,\,\circ\,\,\circ\,\,\circ\,\,\circ\,\,\circ\,\,\circ\,\,\circ\,\,\circ\,\,\circ\,\,\circ\,\,\circ\,\,\circ\,\,\circ\,\,\circ\,\,\circ\,\,\circ\,\,\circ\,\,\circ\,\,\circ\,\,\circ\,\,\circ\,\,\circ\,\,\circ\,\,\circ\,\,\circ\,\,\circ\,\,\circ\,\,\circ\,\,\circ\,\,\circ\,\,\circ\,\,\circ\,\,\bullet\,\,\bullet\,\,\bullet\,\,\bullet\,\,\bullet\,\,\bullet\,\,\bullet\,\,\bullet\,\,\bullet\,\,\bullet\,\,\bullet\,\,\bullet\,\,\bullet\,\,\bullet\,\,\bullet\,\,\bullet\,\,\bullet\,\,\bullet\,\,\bullet\,\,\bullet\,\,\bullet\,\,}\\
^{\circ\,\,\circ\,\,\circ\,\,\circ\,\,\circ\,\,\circ\,\,\circ\,\,\circ\,\,\circ\,\,\circ\,\,\circ\,\,\circ\,\,\circ\,\,\circ\,\,\circ\,\,\circ\,\,\circ\,\,\circ\,\,\circ\,\,\circ\,\,\circ\,\,\circ\,\,\circ\,\,\circ\,\,\circ\,\,\circ\,\,\circ\,\,\circ\,\,\circ\,\,\circ\,\,\circ\,\,\circ\,\,\circ\,\,\circ\,\,\circ\,\,\circ\,\,\bullet\,\,\bullet\,\,\bullet\,\,\bullet\,\,\bullet\,\,\bullet\,\,\bullet\,\,\bullet\,\,\bullet\,\,\bullet\,\,\bullet\,\,\bullet\,\,\bullet\,\,\bullet\,\,\bullet\,\,\bullet\,\,\bullet\,\,\bullet\,\,\bullet\,\,\bullet\,\,}
_{\circ\,\,\circ\,\,\circ\,\,\circ\,\,\circ\,\,\circ\,\,\circ\,\,\circ\,\,\circ\,\,\circ\,\,\circ\,\,\circ\,\,\circ\,\,\circ\,\,\circ\,\,\circ\,\,\circ\,\,\circ\,\,\circ\,\,\circ\,\,\circ\,\,\circ\,\,\circ\,\,\circ\,\,\circ\,\,\circ\,\,\circ\,\,\circ\,\,\circ\,\,\circ\,\,\circ\,\,\circ\,\,\circ\,\,\circ\,\,\circ\,\,\circ\,\,\circ\,\,\bullet\,\,\circ\,\,\bullet\,\,\bullet\,\,\circ\,\,\bullet\,\,\bullet\,\,\bullet\,\,\circ\,\,\bullet\,\,\bullet\,\,\bullet\,\,\bullet\,\,\circ\,\,\bullet\,\,\bullet\,\,\bullet\,\,\bullet\,\,\bullet\,\,}\\
^{\circ\,\,\circ\,\,\circ\,\,\circ\,\,\circ\,\,\circ\,\,\circ\,\,\circ\,\,\circ\,\,\circ\,\,\circ\,\,\circ\,\,\circ\,\,\circ\,\,\circ\,\,\circ\,\,\circ\,\,\circ\,\,\circ\,\,\circ\,\,\circ\,\,\circ\,\,\circ\,\,\circ\,\,\circ\,\,\circ\,\,\circ\,\,\circ\,\,\circ\,\,\circ\,\,\circ\,\,\circ\,\,\circ\,\,\circ\,\,\circ\,\,\circ\,\,\circ\,\,\circ\,\,\bullet\,\,\bullet\,\,\bullet\,\,\bullet\,\,\bullet\,\,\bullet\,\,\bullet\,\,\bullet\,\,\bullet\,\,\bullet\,\,\bullet\,\,\bullet\,\,\bullet\,\,\bullet\,\,\bullet\,\,\bullet\,\,\bullet\,\,\bullet\,\,}
_{\circ\,\,\circ\,\,\circ\,\,\circ\,\,\circ\,\,\circ\,\,\circ\,\,\circ\,\,\circ\,\,\circ\,\,\circ\,\,\circ\,\,\circ\,\,\circ\,\,\circ\,\,\circ\,\,\circ\,\,\circ\,\,\circ\,\,\circ\,\,\circ\,\,\circ\,\,\circ\,\,\circ\,\,\circ\,\,\circ\,\,\circ\,\,\circ\,\,\circ\,\,\circ\,\,\circ\,\,\circ\,\,\circ\,\,\circ\,\,\circ\,\,\circ\,\,\circ\,\,\circ\,\,\circ\,\,\bullet\,\,\bullet\,\,\circ\,\,\bullet\,\,\bullet\,\,\bullet\,\,\circ\,\,\bullet\,\,\bullet\,\,\bullet\,\,\bullet\,\,\circ\,\,\bullet\,\,\bullet\,\,\bullet\,\,\bullet\,\,\bullet\,\,}\\
^{\circ\,\,\circ\,\,\circ\,\,\circ\,\,\circ\,\,\circ\,\,\circ\,\,\circ\,\,\circ\,\,\circ\,\,\circ\,\,\circ\,\,\circ\,\,\circ\,\,\circ\,\,\circ\,\,\circ\,\,\circ\,\,\circ\,\,\circ\,\,\circ\,\,\circ\,\,\circ\,\,\circ\,\,\circ\,\,\circ\,\,\circ\,\,\circ\,\,\circ\,\,\circ\,\,\circ\,\,\circ\,\,\circ\,\,\circ\,\,\circ\,\,\circ\,\,\circ\,\,\circ\,\,\circ\,\,\circ\,\,\bullet\,\,\circ\,\,\circ\,\,\bullet\,\,\bullet\,\,\circ\,\,\circ\,\,\bullet\,\,\bullet\,\,\bullet\,\,\circ\,\,\circ\,\,\bullet\,\,\bullet\,\,\bullet\,\,\bullet\,\,}
_{\circ\,\,\circ\,\,\circ\,\,\circ\,\,\circ\,\,\circ\,\,\circ\,\,\circ\,\,\circ\,\,\circ\,\,\circ\,\,\circ\,\,\circ\,\,\circ\,\,\circ\,\,\circ\,\,\circ\,\,\circ\,\,\circ\,\,\circ\,\,\circ\,\,\circ\,\,\circ\,\,\circ\,\,\circ\,\,\circ\,\,\circ\,\,\circ\,\,\circ\,\,\circ\,\,\circ\,\,\circ\,\,\circ\,\,\circ\,\,\circ\,\,\circ\,\,\circ\,\,\circ\,\,\circ\,\,\circ\,\,\circ\,\,\bullet\,\,\bullet\,\,\bullet\,\,\bullet\,\,\bullet\,\,\bullet\,\,\bullet\,\,\bullet\,\,\bullet\,\,\bullet\,\,\bullet\,\,\bullet\,\,\bullet\,\,\bullet\,\,\bullet\,\,}\\
^{\circ\,\,\circ\,\,\circ\,\,\circ\,\,\circ\,\,\circ\,\,\circ\,\,\circ\,\,\circ\,\,\circ\,\,\circ\,\,\circ\,\,\circ\,\,\circ\,\,\circ\,\,\circ\,\,\circ\,\,\circ\,\,\circ\,\,\circ\,\,\circ\,\,\circ\,\,\circ\,\,\circ\,\,\circ\,\,\circ\,\,\circ\,\,\circ\,\,\circ\,\,\circ\,\,\circ\,\,\circ\,\,\circ\,\,\circ\,\,\circ\,\,\circ\,\,\circ\,\,\circ\,\,\circ\,\,\circ\,\,\circ\,\,\circ\,\,\bullet\,\,\bullet\,\,\bullet\,\,\circ\,\,\bullet\,\,\bullet\,\,\bullet\,\,\bullet\,\,\circ\,\,\bullet\,\,\bullet\,\,\bullet\,\,\bullet\,\,\bullet\,\,}
_{\circ\,\,\circ\,\,\circ\,\,\circ\,\,\circ\,\,\circ\,\,\circ\,\,\circ\,\,\circ\,\,\circ\,\,\circ\,\,\circ\,\,\circ\,\,\circ\,\,\circ\,\,\circ\,\,\circ\,\,\circ\,\,\circ\,\,\circ\,\,\circ\,\,\circ\,\,\circ\,\,\circ\,\,\circ\,\,\circ\,\,\circ\,\,\circ\,\,\circ\,\,\circ\,\,\circ\,\,\circ\,\,\circ\,\,\circ\,\,\circ\,\,\circ\,\,\circ\,\,\circ\,\,\circ\,\,\circ\,\,\circ\,\,\circ\,\,\circ\,\,\bullet\,\,\bullet\,\,\circ\,\,\circ\,\,\bullet\,\,\bullet\,\,\bullet\,\,\circ\,\,\circ\,\,\bullet\,\,\bullet\,\,\bullet\,\,\bullet\,\,}\\
^{\circ\,\,\circ\,\,\circ\,\,\circ\,\,\circ\,\,\circ\,\,\circ\,\,\circ\,\,\circ\,\,\circ\,\,\circ\,\,\circ\,\,\circ\,\,\circ\,\,\circ\,\,\circ\,\,\circ\,\,\circ\,\,\circ\,\,\circ\,\,\circ\,\,\circ\,\,\circ\,\,\circ\,\,\circ\,\,\circ\,\,\circ\,\,\circ\,\,\circ\,\,\circ\,\,\circ\,\,\circ\,\,\circ\,\,\circ\,\,\circ\,\,\circ\,\,\circ\,\,\circ\,\,\circ\,\,\circ\,\,\circ\,\,\circ\,\,\circ\,\,\circ\,\,\bullet\,\,\circ\,\,\circ\,\,\circ\,\,\bullet\,\,\bullet\,\,\circ\,\,\circ\,\,\circ\,\,\bullet\,\,\bullet\,\,\bullet\,\,}
_{\circ\,\,\circ\,\,\circ\,\,\circ\,\,\circ\,\,\circ\,\,\circ\,\,\circ\,\,\circ\,\,\circ\,\,\circ\,\,\circ\,\,\circ\,\,\circ\,\,\circ\,\,\circ\,\,\circ\,\,\circ\,\,\circ\,\,\circ\,\,\circ\,\,\circ\,\,\circ\,\,\circ\,\,\circ\,\,\circ\,\,\circ\,\,\circ\,\,\circ\,\,\circ\,\,\circ\,\,\circ\,\,\circ\,\,\circ\,\,\circ\,\,\circ\,\,\circ\,\,\circ\,\,\circ\,\,\circ\,\,\circ\,\,\circ\,\,\circ\,\,\circ\,\,\circ\,\,\bullet\,\,\bullet\,\,\bullet\,\,\bullet\,\,\bullet\,\,\bullet\,\,\bullet\,\,\bullet\,\,\bullet\,\,\bullet\,\,\bullet\,\,}\\
^{\circ\,\,\circ\,\,\circ\,\,\circ\,\,\circ\,\,\circ\,\,\circ\,\,\circ\,\,\circ\,\,\circ\,\,\circ\,\,\circ\,\,\circ\,\,\circ\,\,\circ\,\,\circ\,\,\circ\,\,\circ\,\,\circ\,\,\circ\,\,\circ\,\,\circ\,\,\circ\,\,\circ\,\,\circ\,\,\circ\,\,\circ\,\,\circ\,\,\circ\,\,\circ\,\,\circ\,\,\circ\,\,\circ\,\,\circ\,\,\circ\,\,\circ\,\,\circ\,\,\circ\,\,\circ\,\,\circ\,\,\circ\,\,\circ\,\,\circ\,\,\circ\,\,\circ\,\,\circ\,\,\bullet\,\,\bullet\,\,\bullet\,\,\bullet\,\,\circ\,\,\bullet\,\,\bullet\,\,\bullet\,\,\bullet\,\,\bullet\,\,}
_{\circ\,\,\circ\,\,\circ\,\,\circ\,\,\circ\,\,\circ\,\,\circ\,\,\circ\,\,\circ\,\,\circ\,\,\circ\,\,\circ\,\,\circ\,\,\circ\,\,\circ\,\,\circ\,\,\circ\,\,\circ\,\,\circ\,\,\circ\,\,\circ\,\,\circ\,\,\circ\,\,\circ\,\,\circ\,\,\circ\,\,\circ\,\,\circ\,\,\circ\,\,\circ\,\,\circ\,\,\circ\,\,\circ\,\,\circ\,\,\circ\,\,\circ\,\,\circ\,\,\circ\,\,\circ\,\,\circ\,\,\circ\,\,\circ\,\,\circ\,\,\circ\,\,\circ\,\,\circ\,\,\circ\,\,\bullet\,\,\bullet\,\,\bullet\,\,\circ\,\,\circ\,\,\bullet\,\,\bullet\,\,\bullet\,\,\bullet\,\,}\\
^{\circ\,\,\circ\,\,\circ\,\,\circ\,\,\circ\,\,\circ\,\,\circ\,\,\circ\,\,\circ\,\,\circ\,\,\circ\,\,\circ\,\,\circ\,\,\circ\,\,\circ\,\,\circ\,\,\circ\,\,\circ\,\,\circ\,\,\circ\,\,\circ\,\,\circ\,\,\circ\,\,\circ\,\,\circ\,\,\circ\,\,\circ\,\,\circ\,\,\circ\,\,\circ\,\,\circ\,\,\circ\,\,\circ\,\,\circ\,\,\circ\,\,\circ\,\,\circ\,\,\circ\,\,\circ\,\,\circ\,\,\circ\,\,\circ\,\,\circ\,\,\circ\,\,\circ\,\,\circ\,\,\circ\,\,\circ\,\,\bullet\,\,\bullet\,\,\circ\,\,\circ\,\,\circ\,\,\bullet\,\,\bullet\,\,\bullet\,\,}
_{\circ\,\,\circ\,\,\circ\,\,\circ\,\,\circ\,\,\circ\,\,\circ\,\,\circ\,\,\circ\,\,\circ\,\,\circ\,\,\circ\,\,\circ\,\,\circ\,\,\circ\,\,\circ\,\,\circ\,\,\circ\,\,\circ\,\,\circ\,\,\circ\,\,\circ\,\,\circ\,\,\circ\,\,\circ\,\,\circ\,\,\circ\,\,\circ\,\,\circ\,\,\circ\,\,\circ\,\,\circ\,\,\circ\,\,\circ\,\,\circ\,\,\circ\,\,\circ\,\,\circ\,\,\circ\,\,\circ\,\,\circ\,\,\circ\,\,\circ\,\,\circ\,\,\circ\,\,\circ\,\,\circ\,\,\circ\,\,\circ\,\,\bullet\,\,\circ\,\,\circ\,\,\circ\,\,\circ\,\,\bullet\,\,\bullet\,\,}\\
^{\circ\,\,\circ\,\,\circ\,\,\circ\,\,\circ\,\,\circ\,\,\circ\,\,\circ\,\,\circ\,\,\circ\,\,\circ\,\,\circ\,\,\circ\,\,\circ\,\,\circ\,\,\circ\,\,\circ\,\,\circ\,\,\circ\,\,\circ\,\,\circ\,\,\circ\,\,\circ\,\,\circ\,\,\circ\,\,\circ\,\,\circ\,\,\circ\,\,\circ\,\,\circ\,\,\circ\,\,\circ\,\,\circ\,\,\circ\,\,\circ\,\,\circ\,\,\circ\,\,\circ\,\,\circ\,\,\circ\,\,\circ\,\,\circ\,\,\circ\,\,\circ\,\,\circ\,\,\circ\,\,\circ\,\,\circ\,\,\circ\,\,\circ\,\,\bullet\,\,\bullet\,\,\bullet\,\,\bullet\,\,\bullet\,\,\bullet\,\,}
_{\circ\,\,\circ\,\,\circ\,\,\circ\,\,\circ\,\,\circ\,\,\circ\,\,\circ\,\,\circ\,\,\circ\,\,\circ\,\,\circ\,\,\circ\,\,\circ\,\,\circ\,\,\circ\,\,\circ\,\,\circ\,\,\circ\,\,\circ\,\,\circ\,\,\circ\,\,\circ\,\,\circ\,\,\circ\,\,\circ\,\,\circ\,\,\circ\,\,\circ\,\,\circ\,\,\circ\,\,\circ\,\,\circ\,\,\circ\,\,\circ\,\,\circ\,\,\circ\,\,\circ\,\,\circ\,\,\circ\,\,\circ\,\,\circ\,\,\circ\,\,\circ\,\,\circ\,\,\circ\,\,\circ\,\,\circ\,\,\circ\,\,\circ\,\,\circ\,\,\bullet\,\,\bullet\,\,\bullet\,\,\bullet\,\,\bullet\,\,}\\
^{\circ\,\,\circ\,\,\circ\,\,\circ\,\,\circ\,\,\circ\,\,\circ\,\,\circ\,\,\circ\,\,\circ\,\,\circ\,\,\circ\,\,\circ\,\,\circ\,\,\circ\,\,\circ\,\,\circ\,\,\circ\,\,\circ\,\,\circ\,\,\circ\,\,\circ\,\,\circ\,\,\circ\,\,\circ\,\,\circ\,\,\circ\,\,\circ\,\,\circ\,\,\circ\,\,\circ\,\,\circ\,\,\circ\,\,\circ\,\,\circ\,\,\circ\,\,\circ\,\,\circ\,\,\circ\,\,\circ\,\,\circ\,\,\circ\,\,\circ\,\,\circ\,\,\circ\,\,\circ\,\,\circ\,\,\circ\,\,\circ\,\,\circ\,\,\circ\,\,\circ\,\,\bullet\,\,\bullet\,\,\bullet\,\,\bullet\,\,}
_{\circ\,\,\circ\,\,\circ\,\,\circ\,\,\circ\,\,\circ\,\,\circ\,\,\circ\,\,\circ\,\,\circ\,\,\circ\,\,\circ\,\,\circ\,\,\circ\,\,\circ\,\,\circ\,\,\circ\,\,\circ\,\,\circ\,\,\circ\,\,\circ\,\,\circ\,\,\circ\,\,\circ\,\,\circ\,\,\circ\,\,\circ\,\,\circ\,\,\circ\,\,\circ\,\,\circ\,\,\circ\,\,\circ\,\,\circ\,\,\circ\,\,\circ\,\,\circ\,\,\circ\,\,\circ\,\,\circ\,\,\circ\,\,\circ\,\,\circ\,\,\circ\,\,\circ\,\,\circ\,\,\circ\,\,\circ\,\,\circ\,\,\circ\,\,\circ\,\,\circ\,\,\circ\,\,\bullet\,\,\bullet\,\,\bullet\,\,}\\
^{\circ\,\,\circ\,\,\circ\,\,\circ\,\,\circ\,\,\circ\,\,\circ\,\,\circ\,\,\circ\,\,\circ\,\,\circ\,\,\circ\,\,\circ\,\,\circ\,\,\circ\,\,\circ\,\,\circ\,\,\circ\,\,\circ\,\,\circ\,\,\circ\,\,\circ\,\,\circ\,\,\circ\,\,\circ\,\,\circ\,\,\circ\,\,\circ\,\,\circ\,\,\circ\,\,\circ\,\,\circ\,\,\circ\,\,\circ\,\,\circ\,\,\circ\,\,\circ\,\,\circ\,\,\circ\,\,\circ\,\,\circ\,\,\circ\,\,\circ\,\,\circ\,\,\circ\,\,\circ\,\,\circ\,\,\circ\,\,\circ\,\,\circ\,\,\circ\,\,\circ\,\,\circ\,\,\circ\,\,\bullet\,\,\bullet\,\,}
_{\circ\,\,\circ\,\,\circ\,\,\circ\,\,\circ\,\,\circ\,\,\circ\,\,\circ\,\,\circ\,\,\circ\,\,\circ\,\,\circ\,\,\circ\,\,\circ\,\,\circ\,\,\circ\,\,\circ\,\,\circ\,\,\circ\,\,\circ\,\,\circ\,\,\circ\,\,\circ\,\,\circ\,\,\circ\,\,\circ\,\,\circ\,\,\circ\,\,\circ\,\,\circ\,\,\circ\,\,\circ\,\,\circ\,\,\circ\,\,\circ\,\,\circ\,\,\circ\,\,\circ\,\,\circ\,\,\circ\,\,\circ\,\,\circ\,\,\circ\,\,\circ\,\,\circ\,\,\circ\,\,\circ\,\,\circ\,\,\circ\,\,\circ\,\,\circ\,\,\circ\,\,\circ\,\,\circ\,\,\circ\,\,\bullet\,\,}\\
\end{bmatrix}
\end{array}$$}

 We say that these sequences as in Table VI are  $q${\it-sequences}, ($q=1,2,3$).
 For $T_6^{11}$, Table VII of $q$-sequences ($q=1,2,3,4$) is obtained.
Let $\sigma_p^q={p+q-1\choose p-1}={p+q-1\choose q}$, for $p,q>0$. These numbers $\sigma_p^q$ form the Pascal's triangle $\Delta$, extended from its portion shown on the left of Table VIII, with actual values on the right enclosure of the table.

The $q$-sequence $s_p^q=\sigma_1^q,\sigma_2^q,\ldots,\sigma_p^q$, ($0<p\le q$), is shown diagonally in $\Delta$, from the upper-left border in the lower-right direction. Its sum (i.e. the sum of its terms) $S(s_p^q)=\sigma_1^q+\sigma_2^q+\ldots+\sigma_p^q=\sigma_p^{q+1}$ is shown in $\Delta$ just down to the left from the last term $\sigma_p^q$ of the sum, because $S(s_p^q)=\sum_{\ell=1}^p{q+\ell\choose\ell-1}={q+p+1\choose p-1}=\sigma_p^{q+1}$.
Thus, the disposition of the numbers $\sigma_p^q$ in the $q$-sequences of a grid ${\cal G}_i^k$ uses the successive sequences $s_1^q,s_2^q,\ldots,s_{i+1}^q$ ($q=0,\ldots,k-i-1$) in the form suggested in Tables VI and VII for ${\cal G}_4^8$ and ${\cal G}_6^{11}$, respectively.

A similar table for the general case ${\cal G}_i^k$ has rows headed by the indices $q=1,2,\ldots,k-i-2,k-i-1$, with the last row (for $q=k-i-1$) just containing the only composing term of the sequence $s_{i+1}^{k-i-1}$, namely the sum $\sigma_{i+1}^{k-i-1}$ of the penultimate row (for $j=k-i-2$) formed by the sequence  $s_{i+1}^{k-i-1}=\sigma_1^{k-i-2},\ldots,\sigma_{i+1}^{k-i-2}$. 

In turn, each term $\sigma_p^{k-i-2}$ of $s_{i+1}^{k-i-1}$ is the sum of the terms of the corresponding sequence $s_p^{k-i-2}$ present in the previous level, and so on.
Table IX lists the values of $\sigma_p^q$ for $p,q=1,\ldots,9$.

\begin{figure}[htp]
\hspace*{1.1cm}%
\includegraphics[scale=0.373]{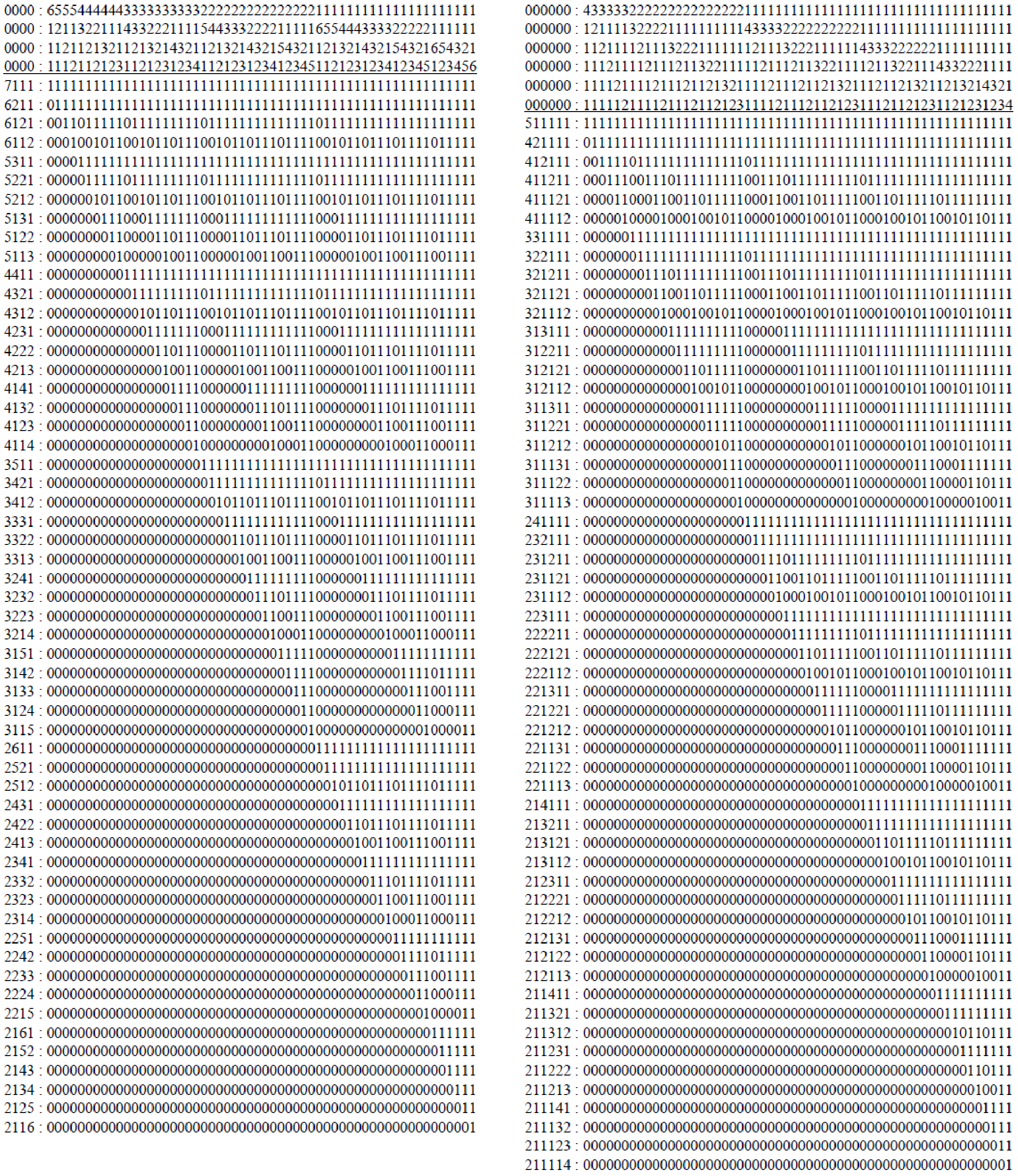}
\caption{$T_4^9$ and $T_6^9$ (row and column headers as in Fig.1)}
\end{figure}

The light-gray components of ${\cal G}_i^k$ (i.e. the maximal unions of side-sharing light-gray single squares) are {\it descending $\ell$-staircases} (i.e. corresponding to the union of all the supports of strictly lower triangular $(1+\ell)\times(1+\ell)$-matrices). Each such staircase has a unique lower-left vertex, namely the intersection of two lines as in items (a)-(b), Remark~\ref{rem}. The number of light-gray single squares in the lower horizontal (or left vertical) side of each such staircase in ${\cal G}_i^k$
 is shown in Table X, where the headers of rows and columns are those indicated in the items (a)-(b), Remark~\ref{rem}.

Tensors $T_4^9$ and $T_5^9$ are represented in Fig.~2, with row and column headers as in Fig.~1. Similarly for tensor $T_5^9$ on the left of Fig.~3. Table XI represents ${\mathcal G}_4^9$.
We state our findings above on distribution of ordered trees in tensors of zipper merging pairs of compositions of contiguous integers in the following statement.

\begin{theorem}\label{teo} The ${k-1\choose i-1}\times{k-1\choose i-1}$-grid ${\cal G}_i^k$
 is partitioned into $q$-blocks, for $q=1,\ldots,i-1$, where each $q$-block is: {\bf(a)} a sub-grid of ${\cal G}_i^k$ with the horizontal and vertical headers of its upper-left (resp. lower-right) corners having (resp. just preceding) the first appearance of an integer $\alpha$
 (resp. $\alpha-1$) in its $(k-q-1)$-entry; {\bf(b)} intersection of an horizontal $q$-strip and a vertical $q$-strip of ${\cal G}_i^k$. The heights and widths of such $q$-strips form $q$-sequences
$s_p^q=\sigma_1^q,\sigma_2^q,\ldots,\sigma_p^q$ ($0<p\le q$)  shown diagonally in $\Delta$ in lower-right direction starting at its upper-left border. Each sum $S(s_p^q)=\sigma_p^{q+1}$ is located in $\Delta$ in the position just down to the left from the last summand $\sigma_p^q$.
The null entries of $T_i^k$ form the union $U_i^k$ of descending staircases $St_B$ in ${\cal G}_i^k$, one per $q$-block $B$, each $St_B$ delimited on the left and below by segments bordering $B$ that intersect at its lower-left corner and have common length $\ell(B)=\eta(B)-1$, with the height $\eta(B)$ of $B$ not larger than its width.
Moreover, the union in
$U_i^k$ can be taken to be disjoint by just including in it those $St_B$'s not contained into larger $q'$-blocks. In fact, each pair of such blocks is either disjoint or one of the two is
strictly contained into the other.
\end{theorem}

The following corollary states that a relation between the diverse tensors $T_i^k$ for fixed $k$ is a direct consequence of Theorem~\ref{teo}.
Two $\ell\times\ell$-matrices $\{t_{p,q}|1\le p\le\ell;1\le q\le\ell\}$ and $\{s_{p,q}|1\le p\le\ell;1\le q\le\ell\}$ will be said to be {\it anti-transpose} whenever
$t_{p,q}=s_{\ell+1-q,\ell+1-p}$, for $1\le p\le\ell$ and $1\le q\le\ell$. A square is said to be {\it anti-symmetric} if it is anti-transpose of itself.

\begin{corollary}\label{theend}
For every $k>1$ it holds that $T_\ell^k$ and $T_{k+1-\ell}^k$ are anti-transpose. In particular, if $\ell=\frac{k+1}{2}$ is even, then $T_\ell^k$ is anti-symmetric.
\end{corollary}

\begin{figure}[htp]
\includegraphics[scale=0.41]{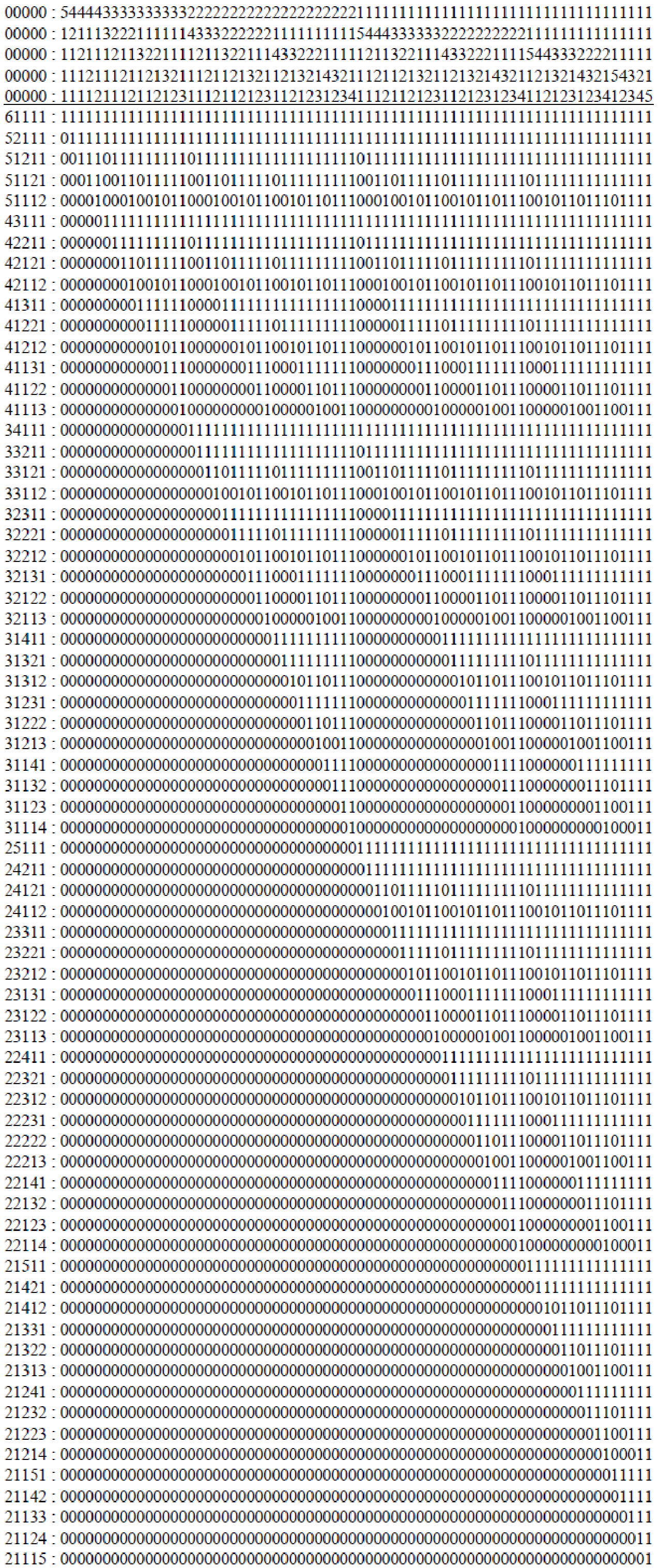}%
\includegraphics[scale=1.73]{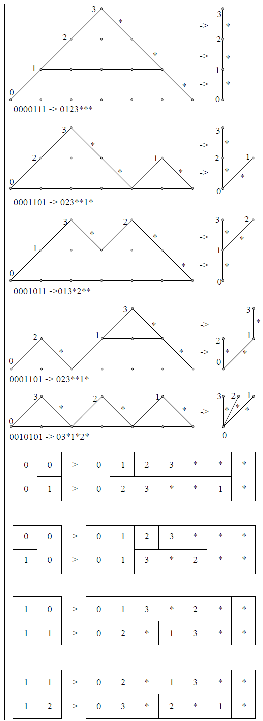}
\caption{$T_5^9$, nested castling via rgs's for $k=3$, and horizontal collapses $\phi_7\rightarrow{\mathcal T}(t)$}
\end{figure}

\end{document}